\documentclass[review,sort&compress,12pt]{elsarticle}
\usepackage{amssymb}
\usepackage{amsmath}
\usepackage{booktabs}
\usepackage{longtable}
\usepackage{graphicx}
\usepackage{float}

\def\listofsymbolsname{List of symbols}
\def\listofsymbols#1{\gdef\@listofsymbols{#1}}

\patchcmd{\maketitle}{\ifx\@elsart#1\endlist}{\ifx\@elsart#1\endlist\par\vskip1em\noindent\textbf{\listofsymbolsname}\par\noindent\@listofsymbols\par}{}{}

\journal{Computer Methods in Applied Mechanics and Engineering}

\begin{document}

\begin{frontmatter}

\title{An Efficient  Explicit-Implicit Adaptive Method for Peridynamic Modelling of Quasi-Static Fracture Formation and Evolution}

\author[CAS]{Shiwei Hu}
\author[CAS,UCAS]{Tianbai Xiao}
\author[CAS,UCAS]{Mingshuo Han}
\author[CAS,UCAS]{Zuoxu Li}
\author[Strathclyde]{Erkan Oterkus}
\author[Strathclyde]{Selda Oterkus}
\author[CAS,UCAS]{Yonghao Zhang\corref{cor1}}

\cortext[cor1]{Corresponding author: yonghao.zhang@imech.ac.cn}

\affiliation[CAS]{organization={Centre for Interdisciplinary Research in Fluids, Institute of Mechanics, Chinese Academy of Sciences},
            city={Beijing},
            postcode={100190}, 
            country={China}}            

\affiliation[UCAS]{organization={School of Engineering Science, University of Chinese Academy of Sciences},
            city={Beijing},
            postcode={100049}, 
            country={China}} 

\affiliation[Strathclyde]{organization={Department of Naval Architecture, Ocean, and Marine Engineering, University of Strathclyde},
            city={Glasgow},
            postcode={G1 1QN}, 
            country={United Kingdom}}

\begin{abstract}
Understanding the quasi-static fracture  formation and evolution is essential for assessing the mechanical properties and structural load-bearing capacity of materials.
Peridynamics (PD) provides an effective computational method to depict fracture mechanics.
The explicit adaptive dynamic relaxation (ADR) method and the implicit methods are two mainstream PD approaches to simulate evolution of quasi-static fractures.
However, no comprehensive and quantitative studies have been reported to compare their accuracy and efficiency.
In this work, we first develop an implicit method for bond-based peridynamics (BBPD) based on the full nonlinear equilibrium equation and the degenerate form of the bond failure function, where the Jacobian matrices are derived using the Newton-Raphson (NR) scheme. Subsequently, we analyze the solvability of the implicit BBPD scheme.
Second, a consistent and comprehensive comparison of  accuracy and efficiency of the explicit ADR and implicit methods is conducted, which reveals computational  efficiency of the implicit methods and their limitations in accurately describing crack formation. Finally, by utilizing the unique advantage of both methods, we develop an adaptive explicit-implicit method and propose a switching criterion to deploy appropriate scheme accordingly.
Four typical quasi-static problems are employed as the numerical experiments, which show the acceleration ratios of the current method range from 6.4 to 141.7 when compared to the explicit ADR.
Therefore, the explicit-implicit adaptive method provides a powerful method to simulate quasi-static fracture formation and evolution.
\end{abstract}

\begin{graphicalabstract}
\includegraphics[clip, width=0.99\textwidth]{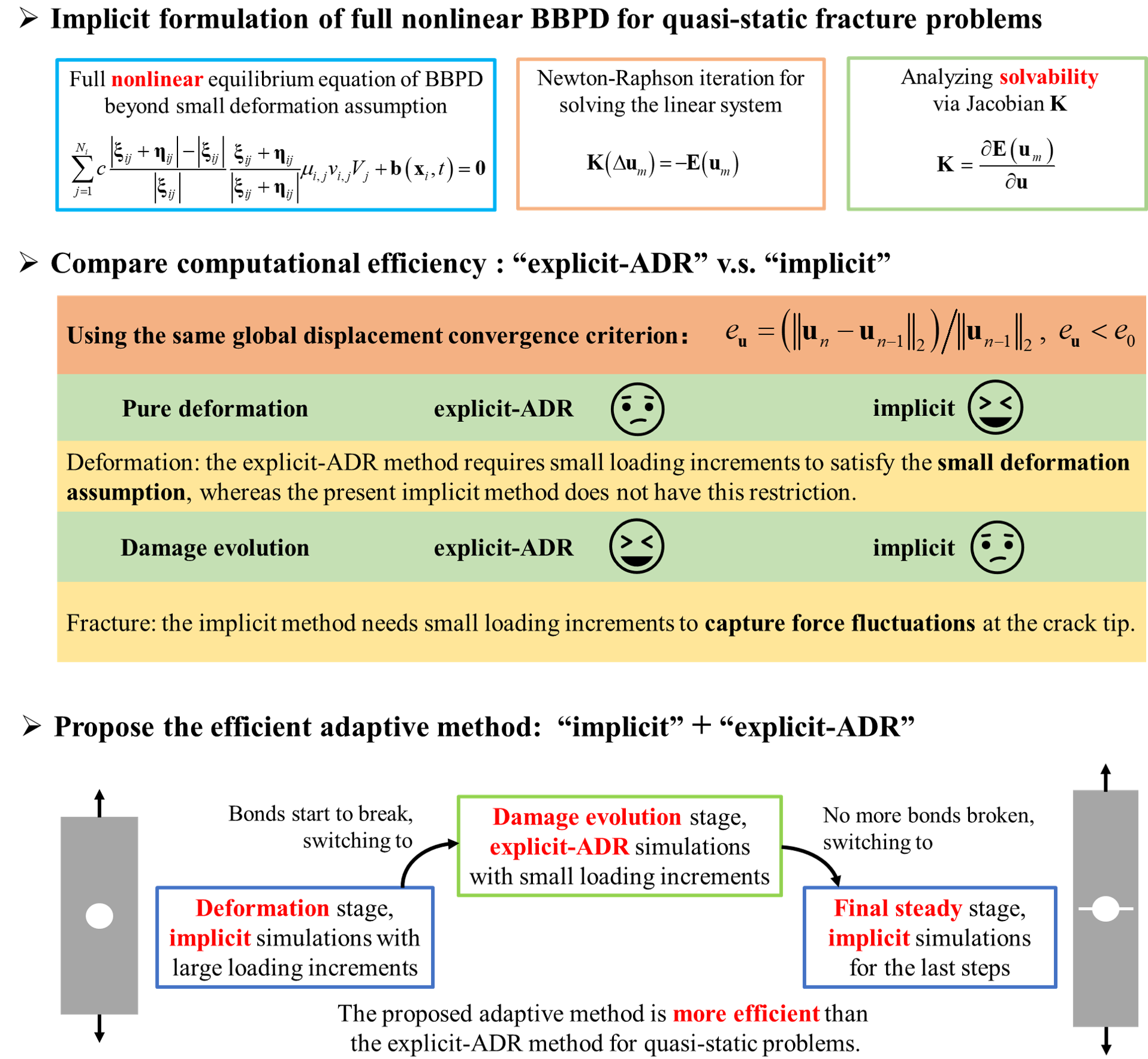}
\end{graphicalabstract}

\begin{highlights}
\item A novel implicit method for the full nonlinear form of bond-based peridynamic model for quasi-static problems
\item A comprehensive comparison of the efficiency and accuracy of explicit adaptive dynamic relaxation and implicit methods
\item An efficient explicit-implicit adaptive method tailored to describe  quasi-static fracture formation and evolution
\end{highlights}


\end{frontmatter}

\newpage

\section*{Nomenclature}
\begin{longtable}{l p{10cm}}
\toprule
\textbf{Symbols} & \textbf{Physical quantities, units} \\
\midrule
\endhead

\bottomrule
\endfoot
$ \mathbf{x} $ & Position vector before deformation, unit: m \\
$ \mathbf{y} $ & Position vector after deformation, unit: m \\
$\mathbf{u}, \ddot{\mathbf{u}}$ & Displacement vector, unit: $m$, acceleration vector, unit: $m/s^2$\\
$\mathbf{\xi}, \mathbf{\eta}$ & Relative position vector, relative displacement vector, unit : $m$\\
$\mathbf{f}(\mathbf{u}^\prime - \mathbf{u},\mathbf{x}^\prime - \mathbf{x},t)$ & Interaction force density vector between the particles at $\mathbf{x'}$ and $\mathbf{x}$ at the time $t$, unit: $N/m^3$ \\
$\mathbf{F}(\mathbf{\xi},\mathbf{u})$ & Set of the total interaction force density vectors acting on each PD point in the entire material, unit: $N/m^3$ \\
$\mathbf{b}$ & Volumetric body force density vector, unit: $N/m^3$\\
$c$ & Bond constant, unit: $Pa/m^4$\\
$s$ & Bond stretch\\
$\rho$ & Density, unit: $Kg/m^3$\\
$E$ & Young’s modulus, unit: $Pa$\\
$N$ & Total number of particles in the computed materials\\
$N_i$ & Number of particles within the horizon $H_i$\\
$\Delta\mathbf{u}_m $ & Displacement increment vector, unit: $m$\\
$\mathbf{K} $ & Jacobian matrix \\
$e$ & Relative error tolerance of the NR procedure\\
$e_0 $ & Relative error tolerance of the  whole-field displacement\\
$S$ & Sparsity index\\
$\Delta x $ & Spatial discretization size, unit: $m$\\
$ T_s$ & Degradation function\\
$s_m $ & Minimum bond stretch for bond interaction degradation\\
$s_c $ & Critical bond stretch\\
$l, w, h $ & Length, width, height, unit: $m$\\
$\Delta u_i $ & Displacement loading increments for the implicit scheme, unit: $m$\\
$\Delta u_e $ & Displacement loading increments for the explicit scheme, unit: $m$\\
$u_{total} $ & Total displacement loading, unit: $m$\\
$N_i $ & Total implicit loading steps to reach the total displacement\\
$N_e $ & Total explicit loading steps to reach the total displacement\\
$r_a $ & Acceleration ratio of the adaptive method in comparison to the explicit-ADR method\\
$r_n $ & Ratio of the duration of the deformation stage to the total simulation duration\\
\end{longtable}

\section{Introduction}
\label{sec1}
Peridynamics (PD) offers an alternative method for modelling and simulation of fracture behavior  \cite{RN1,RN2,RN3,RN4,RN5,RN6,RN7}, which has made rapid progress in the recent years.
Unlike classical continuum mechanics (CCM) that relies on differential operators, PD models the evolution of materials using non-local interactions between the material points within a finite distance.
Thus, it naturally mitigates the well-posedness issue of the derivatives appearing at cracks and voids.
Both continuous and discontinuous fields can be solved within the same framework, where cumbersome remedial techniques are not required \cite{RN9,RN10,RN11,RN12,RN13,RN14,RN15,RN16,RN17,RNzy1,RNzy2}.
Based on modelling of non-local interactions, PD can be divided into two categories, i.e., bond-based PD (BBPD) and state-based PD (SBPD) methods.
BBPD models the interaction between material points as the force of equal magnitude pointing to the opposite directions along the interaction direction  \cite{RN10}.
Despite the limitation of fixed Poisson's ratio, BBPD models are widely used in the modelling of brittle materials due to the favorable numerical stability, succinct mathematical form and low computational cost, and effectiveness in addressing multi-physical coupling problems \cite{RN18,RN19,RN20,RN21,RN22,RN23,RN24,RN25,RN26,RN27,RN28,RN29,RN30,RN31,RN32,RN33,RN56}. Therefore, this study will focus on the BBPD, and may be extended to SBPD in the future.

Quasi-static problems are vital in the study of engineering structural deformation, subsurface movement and deformation, material performance testing, and loading assessment \cite{RN34,RN35,RN36,RN37,RN38,RN39,RN40,RN41}. 
Although PD has been successfully applied to simulate dynamic fracture propagation \cite{RN36,RN37,RN39,RN42,RN43,RN44}, obtaining steady-state solutions under quasi-static conditions remains challenging \cite{RN34}. 
This difficulty arises due to the fact that the small time step employed in the  explicit time integrators leads to  overwhelmingly-costly computation. 
To address this issue, Kilic et al. \cite{RN34} extended the adaptive dynamic relaxation (ADR) method to PD simulation, where the damping coefficient in dynamic relaxation was estimated using Rayleigh's quotient to effectively dampen the system from higher frequency modes to lower frequency modes. This explicit-ADR method can achieve steady-state solutions under quasi-static conditions and accurately capture fracture propagation. However, small loading rates are necessary to avoid nonphysical force fluctuations caused by dynamic effects. \par

On the other hand, the implicit PD schemes have been proposed for static or quasi-static problems \cite{RN35,RN45,RN46,RN47,RN48,RN49,RN50,RN51}, which allow larger loading steps, faster convergence, and guaranteed equilibrium conditions. Among others, Breitenfeld et al. \cite{RN51} presented the development of a static implementation for the non-ordinary state-based PD (NOSBPD) formulation, focusing on small-strain linearly elastic problems. Ni et al. \cite{RN48} proposed an implicit form of the finite element method (FEM)-BBPD coupled model to address static fracture problems and introduced three methods to compute crack propagation, i.e., allowing only one bond to break per step, no control of bond breakage, and limiting the maximum number of bonds that can break at each step. 
It turned out that the third method significantly improved the accuracy of numerical solution. The linearized equilibrium equation was used in that model under the small displacement assumption \cite{RN1,RN52}. Gu et al. \cite{RN56} also applied such linearized equilibrium equation of BBPD within their framework to study porous quasi-brittle materials . If the linearized PD equation is used, small loading increments are required to ensure that the displacements remain sufficiently small at each step to satisfy the prerequisite for linearization, i.e., small deformation.\par

To remove the small deformation restriction, the full nonlinear form of PD is favored in conjunction with iterative schemes such as the Newton-Raphson (NR).
However, the three mainstream crack computation methods discussed above are not suitable for the NR procedure due to the fact that, when there are cracks, the force field becomes discontinuous, making it difficult to obtain the Jacobian matrix. Hashim et al. \cite{RN46} proposed an implicit NOSBPD framework using the NR scheme, where a degradation function is introduced to describe damage evolution, i.e., the interaction degradation of the bond as the bond stretch increases, which is different from immediate disappearance of the bond influence. Yang et al. \cite{RN45} improved the implicit NOSBPD model by incorporating a stress-based criterion to evaluate bond failure, which overcame the difficulties associated with obtaining fracture properties in the previous critical stretch criterion. While the implicit schemes for NOSBPD have been studied, the detailed 
implementation for BBPD, which is equally crucial, has not  garnered adequate attention yet. Also, since implicit methods aim to reduce computational costs, it is necessary to compare the computational efficiency between the explicit-ADR and implicit schemes for quasi-static problems, which has not been systematically examined.\par

Here, we first develop an implicit method based on the full nonlinear equilibrium equation and the degenerate form of the bond failure function of BBPD for quasi-static problems.
The Jacobian matrix is obtained using the NR scheme, which lays the foundation for the implicit scheme.
We thoroughly compare the computational efficiency and accuracy of the  explicit-ADR method and the newly-developed implicit method.
Leveraging the advantages of both methods, we propose an explicit-implicit adaptive method that strikes a balance between accuracy and efficiency for quasi-static fracture problems. 
The performance of the adaptive method will be elaborated for four quasi-static problems.

Note that although implicit schemes were reported  in the literature \cite{RN51, RN45, RN46, RN48, RN56}, the novel contributions of this work lie in two key aspects. First, the nonlinear form of the BBPD equation is retained to eliminate the small-strain restriction often encountered in BBPD studies \cite{RN1, RN34, RN48, RN49, RN52, RN56}. More importantly, this study presents the first comprehensive comparison of two mainstream methods, i.e., the explicit ADR and implicit schemes, for quasi-static problems. Based on this comparison, 
a new adaptive scheme is proposed that effectively combines the advantages of both approaches. Additionally, the current efficient adaptive scheme is also applicable to SBPD.

The rest of the paper is organized as follows. 
Section 2 introduces BBPD and presents a new  implicit model for quasi-static problems.
Section 3 provides a comprehensive comparison between the explicit-ADR method and the newly-developed implicit method.
Section 4 details a novel explicit-implicit adaptive method tailored for quasi-static fracture problems.
Following this, Section 5 highlights four typical case studies, revealing the computational accuracy and efficiency of the explicit-implicit adaptive method. 
The conclusions and further discussion on the remaining issues are presented in the last section.

\section{Theoretical Model}
\label{sec2}

In this section, we formulate a new  implicit BBPD method for quasi-static problems.
The nonlinear equilibrium equation is analyzed by examining the properties of the Jacobian matrix, and the degradation of bond damage is illustrated.

\subsection{Implicit formulation  of full nonlinear BBPD model}
\label{subsec1}

The equation of motion of BBPD is given by 
\begin{equation}\label{eq1}
    \rho(\mathbf{x})\ddot{\mathbf{u}}(\mathbf{x},t) = \int_H\mathbf{f}(\mathbf{u}^\prime - \mathbf{u},
    \mathbf{x}^\prime - \mathbf{x},t)dH + \mathbf{b}(\mathbf{x},t),
\end{equation}
where $H$ represents the horizon zone of the particle at the position $\mathbf{x}$, $\mathbf{f}(\mathbf{u}^\prime - \mathbf{u},
    \mathbf{x}^\prime - \mathbf{x},t)$ is the force density vector between the particles at $\mathbf{x'}$ and $\mathbf{x}$, both at the time $t$, ($\mathbf{x'} - \mathbf{x}$) represents the initial relative position vector between the particles at $\mathbf{x'}$  and $\mathbf{x}$ in the horizon $H$,  $\mathbf{b}$ is the body force density vector at the time $t$, $\rho(\mathbf{x})$ presents the mass density of the particle at $\mathbf{x}$, and $\ddot{\mathbf{u}}$ is the second-order time derivative of displacement $\mathbf{u}$.\par
\begin{figure}[t]
\centering
\includegraphics[width=0.95\textwidth]{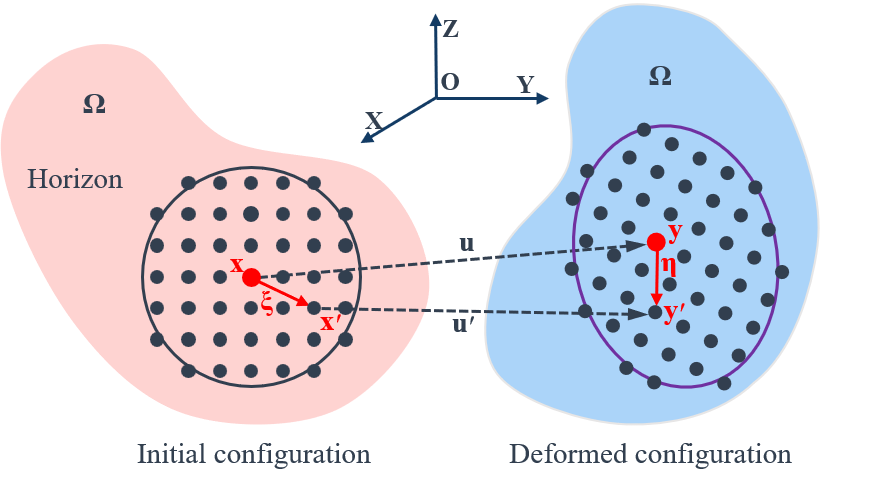}
\caption{Illustration of BBPD and pairwise interaction between PD points. The horizon of the particle at the position $\mathbf{x}$ encompasses all the PD particles interacting with it, and remains unchanged throughout deformation.}\label{fig1}
\end{figure}

The initial relative position vector $\mathbf{\xi}$  and relative displacement vector $\mathbf{\eta}$ are denoted by $\mathbf{\xi}=\mathbf{x'}-\mathbf{x}$ and $\mathbf{\eta}=\mathbf{u'}-\mathbf{u}$, where $\mathbf{u'}$ and $\mathbf{u}$ are, respectively, the displacement of the particles at $\mathbf{x'}$ and $\mathbf{x}$. After the configuration deformation, the new positions of the particles at $\mathbf{x'}$ and $\mathbf{x}$ are represented by $\mathbf{y'}$ and $\mathbf{y}$, where $\mathbf{y'}=\mathbf{x'}+\mathbf{u'}$ and  $\mathbf{y}=\mathbf{x}+\mathbf{u}$. The relative positional relationship between the particles $\mathbf{x}$ and $\mathbf{x'}$ is depicted in Figure \ref{fig1}.
The PD integral-form equation can be discretized for numerical implementation as 
\begin{equation}\label{eq2}
    \rho(\mathbf{x}_i) \ddot{\mathbf{u}}(\mathbf{x}_i,t) = \sum_{j=1}^{N_i} [\mathbf{f}_{i,j}(\mathbf{u}_j - \mathbf{u}_i,
    \mathbf{x}_j - \mathbf{x}_i,t)
    \mu _{i,j}\nu _{i,j} V_j] + \mathbf{b}(\mathbf{x}_i,t),
\end{equation}
where $N_i$ represents the number of particles within the horizon of the particle located at $\mathbf{x}_i$. The subscript $j$ denotes an arbitrary particle positioned at $\mathbf{x}_j$ within the horizon of the particle at $\mathbf{x}_i$. $V_j$ is the  volume of the particle at $\mathbf{x}_j$ within the horizon of the particle at $\mathbf{x}_i$, given by $V_j = (\Delta x)^3$ for a uniform discretization.
Here, $\nu_{i,j}$ and $\mu_{i,j}$ are the volume and surface effect correction factors for the particle at $\mathbf{x}_j$ within the horizon of the particle at $\mathbf{x}_i$, respectively, and $\Delta x$ represents the particle spacing.\par
Based on the BBPD theory, the force density vector of each bond $\mathbf{f}_{i,j}$ satisfies 
\begin{equation}\label{eq3}
\begin{aligned}
    \mathbf{f}_{i,j}(\mathbf{u}_j - \mathbf{u}_i, \mathbf{x}_j - \mathbf{x}_i, t) &= \frac{1}{2} cs \frac{\mathbf{\xi}_{ij} + \mathbf{\eta}_{ij}}{\left| \mathbf{\xi}_{ij} + \mathbf{\eta}_{ij} \right|} \\
    &= \frac{1}{2} c \frac{\left| \mathbf{\xi}_{ij} + \mathbf{\eta}_{ij} \right| - \left| \mathbf{\xi}_{ij} \right|}{\left| \mathbf{\xi}_{ij} \right|} \frac{\mathbf{\xi}_{ij} + \mathbf{\eta}_{ij}}{\left| \mathbf{\xi}_{ij} + \mathbf{\eta}_{ij} \right|} \\
    &= \frac{1}{2} c \frac{\left| \mathbf{y}_j - \mathbf{y}_i \right| - \left| \mathbf{x}_j - \mathbf{x}_i \right|}{\left| \mathbf{x}_j - \mathbf{x}_i \right|} \frac{\mathbf{y}_j - \mathbf{y}_i}{\left| \mathbf{y}_j - \mathbf{y}_i \right|},
\end{aligned}
\end{equation}
where $c$ is bond constant,
\begin{equation*}
    c = \left\{
\begin{array}{rl}
\frac{2E}{\pi\delta^2A} & \text{for\quad 1D,}\\
\frac{9E}{\pi\delta^3h} & \text{for\quad plane stress,}\\
\frac{48E}{5\pi\delta^3h} & \text{for\quad plane strain,}\\
\frac{12E}{\pi\delta^4} & \text{for\quad 3D,}
\end{array} \right.
\end{equation*}
 and $s$ is the bond stretch, which is $s = {{\left( {\left| {{{\bf{\xi }}_{ij}} + {{\bf{\eta }}_{ij}}} \right| - \left| {{{\bf{\xi }}_{ij}}} \right|} \right)} \mathord{\left/
 {\vphantom {{\left( {\left| {{{\bf{\xi }}_{ij}} + {{\bf{\eta }}_{ij}}} \right| - \left| {{{\bf{\xi }}_{ij}}} \right|} \right)} {\left| {{{\bf{\xi }}_{ij}}} \right|}}} \right.
 \kern-\nulldelimiterspace} {\left| {{{\bf{\xi }}_{ij}}} \right|}}$. \par
For the quasi-static deformation, the second-order time derivative of displacement can be ignored, the governing equation of BB-PD is given below,
\begin{equation}
    \sum\limits_{j = 1}^{{N_i}} {c\frac{{\left| {{{\bf{\xi }}_{ij}} + {{\bf{\eta }}_{ij}}} \right| - \left| {{{\bf{\xi }}_{ij}}} \right|}}{{\left| {{{\bf{\xi }}_{ij}}} \right|}}\frac{{{{\bf{\xi }}_{ij}} + {{\bf{\eta }}_{ij}}}}{{\left| {{{\bf{\xi }}_{ij}} + {{\bf{\eta }}_{ij}}} \right|}}} 
    \mu _{i,j}\nu _{i,j} {V_j} 
    + {\bf{b}}\left( {{{\bf{x}}_i},t} \right) = {\bf{0}}.
\end{equation}\par
We propose an implicit algorithm to solve the quasi-static motion equation, thereby avoiding the significant computational cost associated with excessively small time steps required in the explicit time integration. Moreover, the NR scheme is applied to iteratively solve the aforementioned full nonlinear equation.\par
The function $\mathbf{E}(\mathbf{u})$ is defined as $\mathbf{E}(\mathbf{u}) =  \mathbf{F}(\mathbf{\xi}, \mathbf{u}) + \mathbf{b}(\mathbf{x})$, where the function $\mathbf{F}$ represents the set of the total interaction force density vectors acting
on each PD point in the entire material. For the PD point located at $\mathbf{x}_i$, one element of $\mathbf{F}$, the total interaction force density vector at $\mathbf{x}_i$, can be expressed as:
\[\mathbf{F}_i(\mathbf{\xi}, \mathbf{u}) = \sum_{j=1}^{N_i} c \frac{|\mathbf{\xi}_{ij} + \mathbf{\eta}_{ij}| - |\mathbf{\xi}_{ij}|}{|\mathbf{\xi}_{ij}|} \frac{\mathbf{\xi}_{ij} + \mathbf{\eta}_{ij}}{|\mathbf{\xi}_{ij} + \mathbf{\eta}_{ij}|} 
\mu _{i,j}\nu _{i,j} V_j.\]

The system reaches equilibrium when the value of $\mathbf{E}(\mathbf{u})$ approaches zero. \par
\begin{figure}[t]
\centering
\includegraphics[width=0.95\textwidth]{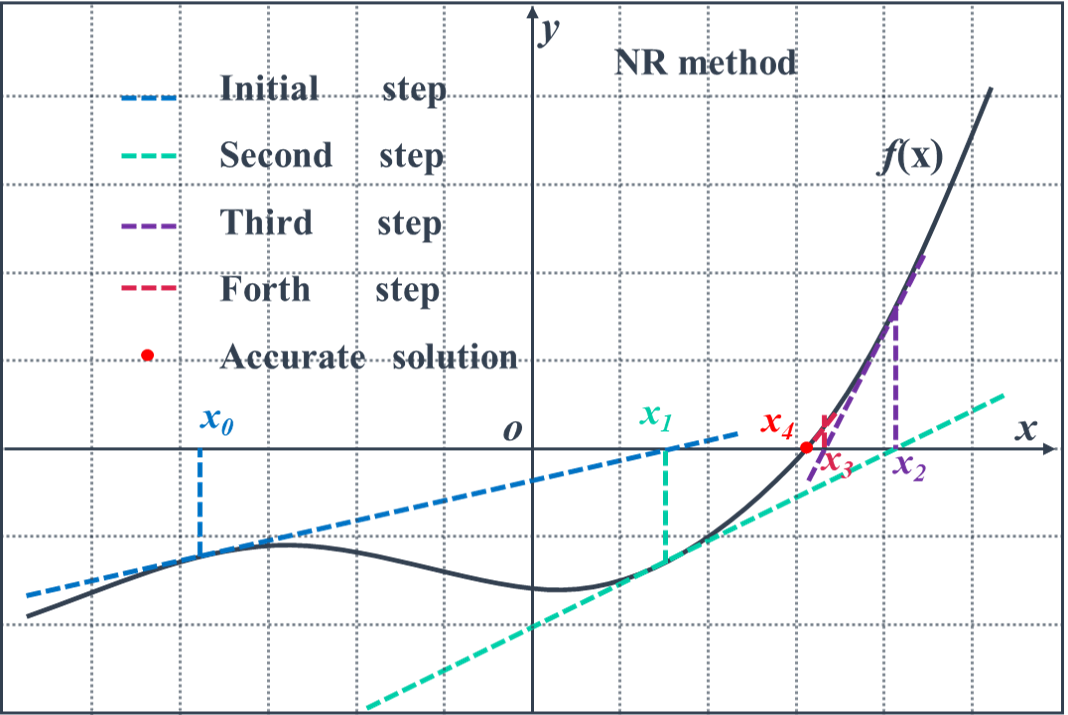}
\caption{Illustration of the NR procedure: the solution of a nonlinear system equals the zero point of the nonlinear function, and the numerical solution obtained through the NR iteration gradually converges toward the accurate solution.}\label{fig2}
\end{figure}

According to the NR iteration scheme shown in figure \ref{fig2}, expand the function $\mathbf{E}(\mathbf{u})$ at $\mathbf{u}_m$ using the Taylor series and retain the linear yields ${\bf{E}}\left( {\mathbf{u}} \right) \approx {\bf{E}}\left( {{{\mathbf{u}}_m}} \right) + \frac{{\partial {\bf{E}}}}{{\partial {\mathbf{u}}}}({\bf{u}}_m)\left( {\mathbf{u}-\mathbf{u}_m} \right)$, where $\mathbf{u}_{m}$ represents the displacement of the current $m$ iteration step. The solution to the nonlinear system $\mathbf{E}(\mathbf{u}) = \mathbf{0}$ can be transformed into solving the following linear equation at each iteration step, 
\begin{equation}\label{eq5}
{{\bf{K}}_m}\left( {\Delta {{\bf{u}}_m}} \right) =  - {\bf{E}}\left( {{{\bf{u}}_m}} \right),{\rm{ }}{{\bf{K}}_m} = \frac{{\partial {\bf{E}}}}{{\partial {\bf{u}}}}(\mathbf{u}_m),
\end{equation}
where $\mathbf{K}_m$ is the Jacobian matrix (or tangent stiffness matrix) at the $m$th iteration, and $\Delta \mathbf{u}_m = \mathbf{u}_{m+1} - \mathbf{u}_m$ is the unknown displacement increment. Once $\Delta \mathbf{u}_m$ is determined, the next displacement $\mathbf{u}_{m+1}$ is updated as $\mathbf{u}_{m+1} = \mathbf{u}_m + \Delta \mathbf{u}_m$. This iterative process continues until $\mathbf{u}_{m+1}$ converges to the required accuracy.

Although ${\bf{K}}_m$ can be determined numerically \cite{RN50}, the calculation involves finite differencing of the whole-field function $\mathbf{E}(\mathbf{u})$ with respect to the displacements of all PD points (as shown in Eq.(\ref{eq5})), resulting in significant computational and memory costs. Furthermore, $\mathbf{K}_m$ changes with each iteration, solving it numerically becomes increasingly expensive.
To reduce the computational cost, the analytical
 expression of $\bf{K}$ can be obtained as,
\begin{equation}\label{eq6}
\bf{K} = \frac{{\partial \bf{E}}}{{\partial {\bf{u}}}} = \frac{{\partial {\bf{F}}\left( {{\bf{\xi }},{\bf{u}}} \right)}}{{\partial {\bf{u}}}} = {\left\lfloor {\begin{array}{*{20}{c}}
{\frac{{\partial {{\bf{F}}_1}\left( {{\bf{\xi }},{\bf{u}}} \right)}}{{\partial {{\bf{u}}_1}}}}&{\frac{{\partial {{\bf{F}}_1}\left( {{\bf{\xi }},{\bf{u}}} \right)}}{{\partial {{\bf{u}}_2}}}}& \cdots &{\frac{{\partial {{\bf{F}}_1}\left( {{\bf{\xi }},{\bf{u}}} \right)}}{{\partial {{\bf{u}}_n}}}}\\
{\frac{{\partial {{\bf{F}}_2}\left( {{\bf{\xi }},{\bf{u}}} \right)}}{{\partial {{\bf{u}}_1}}}}&{\frac{{\partial {{\bf{F}}_2}\left( {{\bf{\xi }},{\bf{u}}} \right)}}{{\partial {{\bf{u}}_2}}}}& \cdots &{\frac{{\partial {{\bf{F}}_2}\left( {{\bf{\xi }},{\bf{u}}} \right)}}{{\partial {{\bf{u}}_N}}}}\\
 \vdots & \vdots & \ddots & \vdots \\
{\frac{{\partial {{\bf{F}}_N}\left( {{\bf{\xi }},{\bf{u}}} \right)}}{{\partial {{\bf{u}}_1}}}}&{\frac{{\partial {{\bf{F}}_N}\left( {{\bf{\xi }},{\bf{u}}} \right)}}{{\partial {{\bf{u}}_2}}}}& \cdots &{\frac{{\partial {{\bf{F}}_N}\left( {{\bf{\xi }},{\bf{u}}} \right)}}{{\partial {{\bf{u}}_N}}}}
\end{array}} \right\rfloor _{N \times N}},
\end{equation}
where $\bf{K}$ represents the Jacobian matrix function at arbitrary time and $N$ represents the total number of particles of the material. For three dimensions, the 3×3 matrix is given by
\begin{equation}\label{eq7}
\frac{{\partial {{\bf{F}}_i}\left( {{\bf{\xi }},{\bf{u}}} \right)}}{{\partial {{\bf{u}}_j}}} = {\left[ {\begin{array}{*{20}{c}}
{\frac{{\partial {F_{ix}}\left( {{\bf{\xi }},{\bf{u}}} \right)}}{{\partial {u_{jx}}}}}&{\frac{{\partial {F_{ix}}\left( {{\bf{\xi }},{\bf{u}}} \right)}}{{\partial {u_{jy}}}}}&{\frac{{\partial {F_{ix}}\left( {{\bf{\xi }},{\bf{u}}} \right)}}{{\partial {u_{jz}}}}}\\
{\frac{{\partial {F_{iy}}\left( {{\bf{\xi }},{\bf{u}}} \right)}}{{\partial {u_{jx}}}}}&{\frac{{\partial {F_{iy}}\left( {{\bf{\xi }},{\bf{u}}} \right)}}{{\partial {u_{jy}}}}}&{\frac{{\partial {F_{iy}}\left( {{\bf{\xi }},{\bf{u}}} \right)}}{{\partial {u_{jz}}}}}\\
{\frac{{\partial {F_{iz}}\left( {{\bf{\xi }},{\bf{u}}} \right)}}{{\partial {u_{jx}}}}}&{\frac{{\partial {F_{iz}}\left( {{\bf{\xi }},{\bf{u}}} \right)}}{{\partial {u_{jy}}}}}&{\frac{{\partial {F_{iz}}\left( {{\bf{\xi }},{\bf{u}}} \right)}}{{\partial {u_{jz}}}}}
\end{array}} \right]_{3 \times 3}},
\end{equation}
where the subscript $x$, $y$, and $z$ represent three mutually orthogonal Cartesian coordinate components, respectively.\par
If we substitute Eq.(\ref{eq7}) into Eq.(\ref{eq6}) for $j \neq i$, and the particle at the position $\mathbf{x}_j$  belongs to the horizon of the particle at the position $\mathbf{x}_i$, i.e., $\mathbf{x}_j \in H_i$, we obtain 
\begin{equation}\label{eq8}
\small
\frac{{\partial {F_{ip}}\left( {{\bf{\xi }},{\bf{u}}} \right)}}{{\partial {u_{jq}}}} = 
c\mu _{i,j}\nu _{i,j} {V_j}
\left[ {\frac{\delta_{pq}}{{\left| {{{\bf{x}}_j} - {{\bf{x}}_i}} \right|}} - \frac{\delta_{pq}}{{\left| {{{\bf{y}}_j} - {{\bf{y}}_i}} \right|}} + \frac{{\left( {{y_{jp}} - {y_{ip}}} \right)\left( {{y_{jq}} - {y_{iq}}} \right)}}{{{{\left| {{{\bf{y}}_j} - {{\bf{y}}_i}} \right|}^3}}}} \right], 
\end{equation}
 where the subscripts $i$ and $j$ represent the coordinate indices ($i, j = 1, 2, ... N$), and $p$ and $q$ denote the indices for the orthogonal component directions of the position vectors ($p, q = x, y, z$). For example, $y_{jp}$ represents the position coordinate component of particle $j$ in the $p$ direction after deformation. $\delta_{pq}=1$ when $p=q$, and $\delta_{pq}=0$ otherwise.\par
If $j \neq i$ and the particle at the position $\mathbf{x}_j$ does not belong to the horizon of the particle at the position $\mathbf{x}_i$,i.e., $\mathbf{x}_j \notin H_i$, we can obtain 
\begin{equation}\label{eq9}
\frac{{\partial {{\bf{F}}_i}\left( {{\bf{\xi }},{\bf{u}}} \right)}}{{\partial {{\bf{u}}_j}}} = {\bf{0}}.
\end{equation}
If $j = i$, similarly, we obtain 
\begin{equation}\label{eq10}
\small
\frac{{\partial {F_{ip}}\left( {{\bf{\xi }},{\bf{u}}} \right)}}{{\partial {u_{iq}}}} = 
c \sum\limits_{\mathbf{x}_k \in H_i} {\left[ {\frac{\delta_{pq}}{{\left| {{{\bf{y}}_k} - {{\bf{y}}_i}} \right|}} - \frac{\delta_{pq}}{{\left| {{{\bf{x}}_k} - {{\bf{x}}_i}} \right|}} - \frac{{\left( {{y_{kp}} - {y_{ip}}} \right)\left( {{y_{kq}} - {y_{iq}}} \right)}}{{{{\left| {{{\bf{y}}_k} - {{\bf{y}}_i}} \right|}^3}}}} \right]\mu _{i,k}\nu _{i,k} {V_k}},
\end{equation}
where the subscript $k$ represents the coordinate indices, and $k = 1, 2, ... N$. 

In summary, the dimension of Jacobian matrix $\bf{K}$ is $mN \times mN$, where $m$ represents the dimension of the system ranging from 1 to 3. Each element of the Jacobian matrix $\bf{K}$ can be calculated by Eqs.(\ref{eq6})-(\ref{eq10}).\par
In the same manner, the right-hand side of  Eq.(\ref{eq5}), which is a $mN \times 1$ vector, can be calculated by 
\begin{equation}\label{eq11}
\begin{aligned}
    - E\left( {{u_{ip}}} \right) = &\ 
    c \sum\limits_{\mathbf{x}_k \in H_i} \left[ \frac{{\left( {{y_{kp}} - {y_{ip}}} \right)}}{{\left| {{{\bf{y}}_k} - {{\bf{y}}_i}} \right|}} - \frac{{\left( {{y_{kp}} - {y_{ip}}} \right)}}{{\left| {{{\bf{x}}_k} - {{\bf{x}}_i}} \right|}} \right]\mu _{i,k}\nu _{i,k} {V_k}  - {b_{ip}}, \\
    & \left(i,k = 1,2 \ldots N;\quad p,q = x,y,z \right).
\end{aligned}
\end{equation}\par
By solving the linear system, $\Delta {\bf{u}}_m$ can be obtained, and the displacement of the next step is $\mathbf{u}_{m+1}=\mathbf{u}_{m}+\Delta\mathbf{u}_m$. The convergence criterion of the NR procedure is set as $\frac{{{{\left\| {{\bf{E}}\left( {{{\bf{u}}_{m + 1}}} \right)} \right\|}_2}}}{{{{\left\| {{\bf{b}}\left( {\bf{\xi }} \right)} \right\|}_2}}} \le e$,  where $e$ stands for the relative error tolerance.\par

\subsection{Properties of Jacobian matrix for implicit BBPD model}
\label{subsec2}
The typical Jacobian matrix of BBPD for quasi-static problems exhibits the following properties.\par
(1) Sparsity:  the number of non-zero elements in $\bf{K}$ is associated with the total number of particles in a horizon. The sparsity index $S$ can be calculated as  $S = {{\sum\limits_{i = 1}^N {\left( {{N_i} + 1} \right)} } \mathord{\left/
 {\vphantom {{\sum\limits_{i = 1}^N {\left( {{N_i} + 1} \right)} } {{N^2}}}} \right.
 \kern-\nulldelimiterspace} {{N^2}}}$,  where ($N_i$+1) represents the total number of particles in the horizon of particle $i$, including itself. Furthermore, not all horizons contain the same number of particles due to the presence of surfaces, so $S < \frac{N_m + 1}{N}$, where $N_m$ is the maximum number of localized particles, which is typically constant in PD simulations. As the horizon size is commonly set to be $3\Delta x$, the value of $N_m$ is typically equal to 3 for 1-dimensional (1D) problems, 28 for 2D problems, and 122 for 3D problems. Therefore, $K$ is sparse as $N \gg N_m$.\par
 (2) Symmetry: if we take the 3D Jacobian matrix $\bf{K}$ as an example, each element of $\mathbf{K}_{3N \times 3N}$ can be expressed as $[$K$_{i,j}]_{p,q}$, where $i$ and $j$ represent the index numbers of all the particles, and $p$ and $q$ represent the axis components of the Cartesian coordinate.\par
 If $j \neq i$, from Eq.(\ref{eq8}), we obtain
\begin{equation*}
\begin{aligned}
    {{{\left[ {{K_{i,j}}} \right]}_{p,q}}} &= 
    c\mu _{i,j}\nu _{i,j} {V_j}
    \left[ \frac{\delta_{pq}}{\left| {{{\bf{x}}_j} - {{\bf{x}}_i}} \right|} - \frac{\delta_{pq}}{\left| {{{\bf{y}}_j} - {{\bf{y}}_i}} \right|} + \frac{\left( {{y_{jp}} - {y_{ip}}} \right)\left( {{y_{jq}} - {y_{iq}}} \right)}{\left| {{{\bf{y}}_j} - {{\bf{y}}_i}} \right|^3} \right], \\
    {{{\left[ {{K_{j,i}}} \right]}_{q,p}}} &=
    c\mu _{i,j}\nu _{i,j} {V_i}
    \left[ \frac{\delta_{pq}}{\left| {{{\bf{x}}_i} - {{\bf{x}}_j}} \right|} - \frac{\delta_{pq}}{\left| {{{\bf{y}}_i} - {{\bf{y}}_j}} \right|} + \frac{\left( {{y_{iq}} - {y_{jq}}} \right)\left( {{y_{ip}} - {y_{jp}}} \right)}{\left| {{{\bf{y}}_i} - {{\bf{y}}_j}} \right|^3} \right].
\end{aligned}
\end{equation*}
In PD simulations, the grids are usually uniform, which mean $V_j = V_i$, so that, $[$K$_{i,j}]_{p,q} = [$K$_{j,i}]_{q,p} $. 
If $j = i$, the elements are on the diagonal of a square matrix, so $[$K$_{i,i}]_{p,q} = [$K$_{i,i}]_{q,p} $.
 Therefore, $\bf{K}$ is a symmetric matrix. \par
(3) Other important property: as for the diagonal elements of the particle index matrix [$K_{i,i}$], the arbitrary Cartesian coordinate components [$K_{i,i}$]$_{p,q}$ can be expressed as Eq.(\ref{eq10}), and we get ${\left[ {{K_{i,i}}} \right]_{p,q}} =  - \sum\limits_{j = 1}^{{N_i}} {{{\left[ {{K_{i,j}}} \right]}_{p,q}}}$. This implies that $\bf{K}$ is a matrix with the row sum equal to 0. Therefore, $\bf{K}$ is singular and requires the introduction of boundary condition to eliminate singularity.\par
\begin{figure}[t]
\centering
\includegraphics[width=0.8\textwidth]{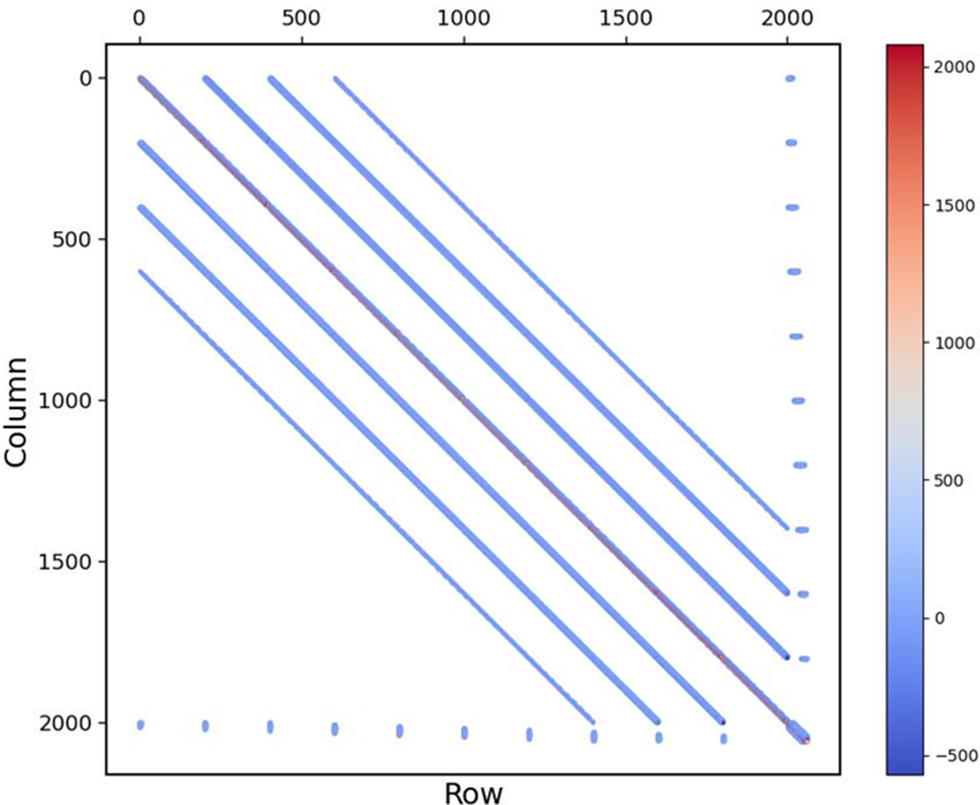}
\caption{Distribution of non-zero elements in the  stiffness matrix shows the unique sparse band-like  characteristics.}\label{fig3}
\end{figure}

To be more illustrative, we compute a specific Jacobian matrix for a particular quasi-static case before applying boundary conditions. This case involves a 2D bar subjected to a transverse loading, which will be discussed in Section \ref{subsec3.1}. Figure \ref{fig3} displays the distribution of non-zero elements in this Jacobian matrix.
As for the positive definiteness of the Jacobian matrix, the theoretical proof is yet to be found while the numerical verification has accomplished.\par 

\subsection{Boundary condition}
\label{subsec2.3}
\begin{figure}[t]
\centering
\includegraphics[width=0.8\textwidth]{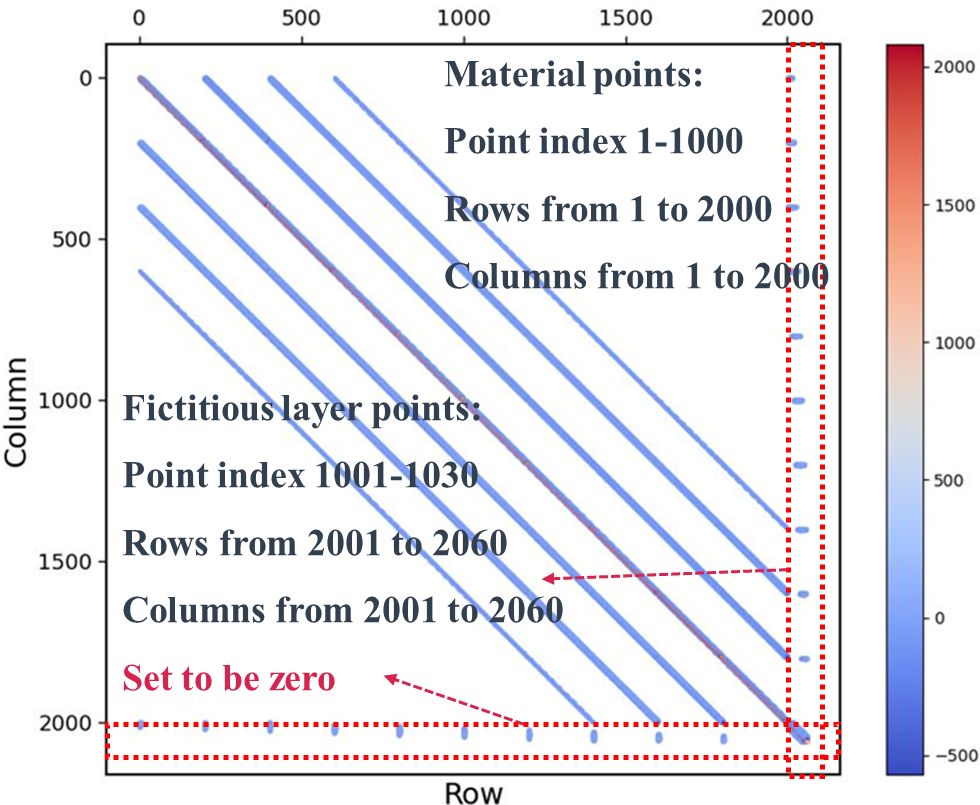}
\caption{Illustration of applying the boundary condition for a 2D bar under a  transverse loading, similar to the scenario depicted in Figure \ref{fig3}. The boundary condition involves fixing the left side of bar.}\label{fig4}
\end{figure}

The Jacobian matrix $\bf{K}$ may become singular due to the presence of rigid body displacement. In order to address this singularity, displacement constraints need to be applied. In PD simulations, a displacement boundary condition is enforced by imposing constraints on the displacement or velocity field within a fictitious material layer along the boundary of a non-zero volume. As the displacement of the materials points in the fictitious layer is known, the rows and columns corresponding to the material points in the fictitious layer should be removed from the Jacobian matrix. Therefore, the row sums are no longer zero, and the singularity of $\bf{K}$ is eliminated,  which is shown in Figure \ref{fig4}.\par

\subsection{Damage evolution}
\label{subsec2.4}
A critical bond stretch is used to describe the damage evolution in the material. When the bond stretch exceeds this critical value, bonds begin to break, cracks start to form, and the interaction forces between the bonds are then set to zero. However, it is important to note that directly setting the forces to zero may adversely affect numerical convergence in the present implicit procedure. The current implicit format is based on the NR scheme, and the differentiation of forces with respect to displacement necessitates the continuity of forces. Hence, the continuous degradation function $T_s$ is utilized to depict the deterioration of bond interaction as the bond stretch increases, which is shown below:
\begin{equation}\label{eq12}
    {T_s} = \left\{ {\begin{array}{r}
{1 \qquad \qquad {\rm{                                                  }}\left( {s \le {s_m}} \right)}\\
{\frac{1}{2}\left[ {1 - \tanh \frac{{\beta \left( {{s_m} + {s_c} - 2s} \right)}}{{{s_m} - {s_c}}}} \right]{\rm{  }}  \left( {{s_m} < s < {s_c}} \right)}\\
{{\rm{0                                                   }}  \qquad \qquad \left( {s \ge {s_c}} \right)},
\end{array}} \right.
\end{equation}
where $s_m$ is the minimum bond stretch at which the degradation of bond interaction begins, $s_c$ is the critical bond stretch, and $\beta$ is a non-negative value that controls the rate of degradation\cite{RN46}.\par
The new total interaction force density vector of the particle at $\mathbf{x}_i$, taking into account the damage, is given by the following equation:
\begin{equation}\label{eq13}
    {{\bf{F}}_i}\left( {{{\bf{\xi }}_{ij}},{{\bf{u}}_{ij}}} \right) = \sum\limits_{j = 1}^{{N_i}} {cs\frac{{{{\bf{\xi }}_{ij}} + {{\bf{\eta }}_{ij}}}}{{\left| {{{\bf{\xi }}_{ij}} + {{\bf{\eta }}_{ij}}} \right|}}} {T_s}\mu _{i,j}\nu _{i,j} {V_j}.
\end{equation}\par
By conducting the NR procedure, the elements of new Jacobian matrix $\mathbf{K}$ can be calculated by
\begin{equation}\label{eq14}
\begin{aligned}
    \frac{\partial F_{ip}({\bf{\xi }},{\bf{u}})}{\partial u_{jq}} &= T_s P_{i,j} A_{ijpq} + P_{i,j} \left( T_s C_{ij} + L_{ij} B_{ij} A_{ijpq} \right) D_{ijpq}, \\
    & \left( i,j = 1,2,\ldots,N,\ i \ne j; \quad p,q = x,y,z \right),
\end{aligned}
\end{equation}
and 
\begin{equation}\label{eq15}
    \begin{aligned}
    \frac{{\partial {F_{ip}}\left( {{\bf{\xi }},{\bf{u}}} \right)}}{{\partial {u_{iq}}}} = &\ - \sum\limits_{\mathbf{x}_j \in H_i} \left[ {T_s}{P_{i,j}}{A_{ijpq}} + {P_{i,j}}\left( {T_s}{C_{ij}} + {L_{ij}}{B_{ij}}{A_{ij}} \right){D_{ijpq}} \right], \\
    &\left(i,j = 1,2,\ldots, N; \quad p,q = x,y,z \right),
\end{aligned}
\end{equation}
where ${P_{i,j}} = c\mu _{i,j}\nu _{i,j} {V_j},{\rm{ }}{A_{ijpq}} = \frac{\delta_{pq}}{{\left| {{{\bf{x}}_j} - {{\bf{x}}_i}} \right|}} - \frac{\delta_{pq}}{{\left| {{{\bf{y}}_j} - {{\bf{y}}_i}} \right|}},{\rm{ }}{B_{ij}} = \frac{1}{{\left| {{{\bf{x}}_j} - {{\bf{x}}_i}} \right|\left| {{{\bf{y}}_j} - {{\bf{y}}_i}} \right|}},{\rm{ }}{C_{ij}} = \frac{1}{{{{\left| {{{\bf{y}}_j} - {{\bf{y}}_i}} \right|}^3}}}$ and ${D_{ijpq}} = \left( {{y_{jp}} - {y_{ip}}} \right)\left( {{y_{jq}} - {y_{iq}}} \right),{\rm{ }}{L_{ij}} = \frac{\beta }{{{s_m} - {s_c}}}\left( {1 - {{\tanh }^2}\frac{{\beta \left( {{s_m} + {s_c} - 2s} \right)}}{{{s_m} - {s_c}}}} \right)$.\par
The properties of Jacobian matrix are not affected by the degradation function $T_s$ and remain unchanged.\par
In summary, by introducing the degradation damage function and applying the boundary condition, the Jacobian  matrix holds its sparse, symmetric, positive definite, and non-singular properties, indicating good solvability of the system.\par

\section{Comparison of explicit-ADR and implicit BBPD method}
\label{sec3}
The implicit scheme is often designed to reduce computational costs of PD for quasi-static problems. With the deploy of BBPD, which represents the simplest form of PD theory, and the implementation of conjugate gradient (CG) method, known as one of the simplest and most efficient iterative methods for solving linear equations, our implicit scheme is a promising method. Therefore, it would be persuasive to evaluate the efficiency of our implicit scheme by comparing 
to the explicit-ADR method. \par
For consistency, the whole-field displacement convergence is maintained for both the explicit-ADR and implicit schemes, as described below:
\begin{equation}\label{eq16}
    {e_{\bf{u}}} = {{\left( {{{\left\| {{{\bf{u}}_n} - {{\bf{u}}_{n - 1}}} \right\|}_2}} \right)} \mathord{\left/
 {\vphantom {{\left( {{{\left\| {{{\bf{u}}_n} - {{\bf{u}}_{n - 1}}} \right\|}_2}} \right)} {{{\left\| {{{\bf{u}}_{n - 1}}} \right\|}_2}}}} \right.
 \kern-\nulldelimiterspace} {{{\left\| {{{\bf{u}}_{n - 1}}} \right\|}_2}}},{\rm{ }}{e_{\bf{u}}} < {e_0},
\end{equation}
where $e_{\bf{u}}$ is the defined relative error of whole-field displacement at the step $n$, $e_0$ is the converge criterion,  $\bf{u}_n$ and $\bf{u}_{n-1}$ are the displacement at the steps ($n$-1) and $n$, respectively. The same material properties, geometric information, discretization length, boundary conditions, initial loading conditions, and displacement convergence criteria are used. Additionally, the degradation function with the same fracture criterion is applied to describe the bond failure in both simulations. All the calculations are conducted using Julia\cite{RN53} programming within the same computational environment, utilizing a single core of the AMD EPYC 7763 64-Core Processor. Subsequently, the results of damage, displacement, and CPU time for computations are compared.\par

\subsection{Deformation of undamaged 2D bar}
\label{subsec3.1}
\begin{figure}[t]
\centering
\includegraphics[width=0.95\textwidth]{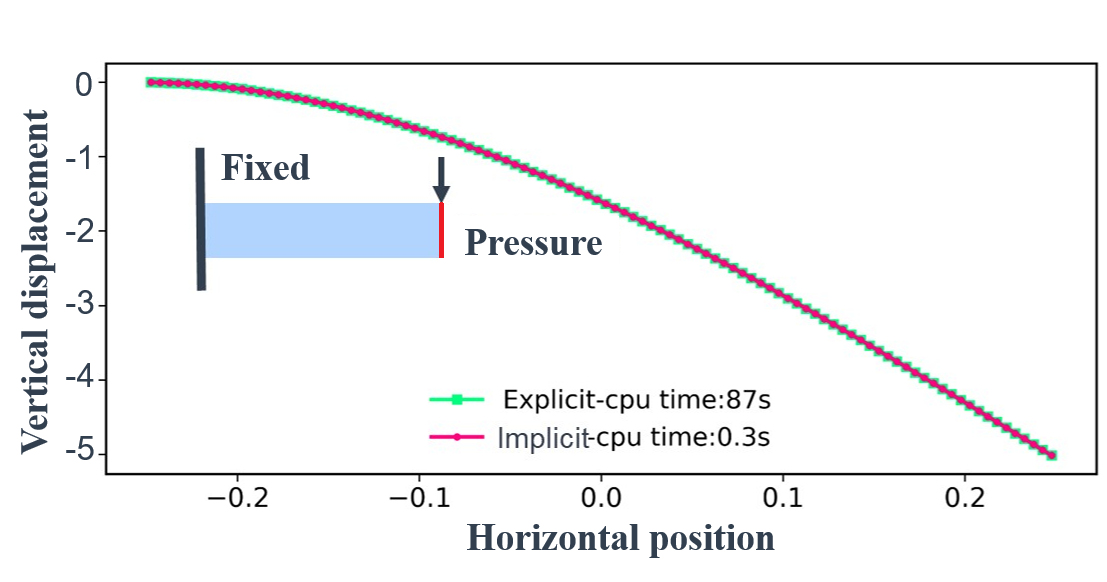}
\caption{This figure illustrates the vertical displacement distribution along the horizontal central line of a 2D bar subjected to transverse loads. The results are obtained from the explicit-ADR and implicit simulations, respectively. The loading position is marked in red.}\label{fig5}
\end{figure}

The 2D bar has the initial length $l$ of 0.5m, and the width $w$ of 0.05m. It is uniformly discretized into 1030 particles using the regular grid 100 × 10, with a layer of fictitious thickness 3$\Delta x$ at the fixed boundary on the left side of the bar. The horizon size is 3.015$\Delta x$. The material is assumed to be isotropic and linear elastic with the Young’s modulus $E$ of 200GPa and the density of 7850 kg/$m^3$. The boundary and loading conditions are depicted in Figure \ref{fig5}. The left end of the bar is fixed, while the right end is subjected to a downward vertical loading of 125 N. The loading boundary is located on the outermost layer of the right end of the bar, with a boundary layer volume of 10×1×$(\Delta x)^3$. The loading is applied from the beginning as a body force density vector to the PD points within the loading boundary layer, with a magnitude of 1×$10^8$ N/$m^3$.\par

Both the implicit and explicit schemes adhere to the same convergence criterion as described in Eq.(\ref{eq16}), and $e_0$ is set to be 1×$10^{-9}$. The implicit simulation converges in only one step, while the explicit one requires approximately 20000 steps. This is because small loading increment is required to satisfy the small displacement assumption of the explicit-ADR method\cite{RN10}, leading to a large number of steps, whereas it is not required 
for the present implicit method. \par
No damage is allowed in this scenario. As depicted in Figure \ref{fig5}, the results from the implicit simulation are identical to those obtained from the explicit simulations. But the implicit scheme is approximately 290 times faster than the explicit one.\par
\subsection{Fracture of 2D plate with a hole}
\label{subsec3.2}
\begin{figure}[t]
\centering
\includegraphics[width=0.99\textwidth]{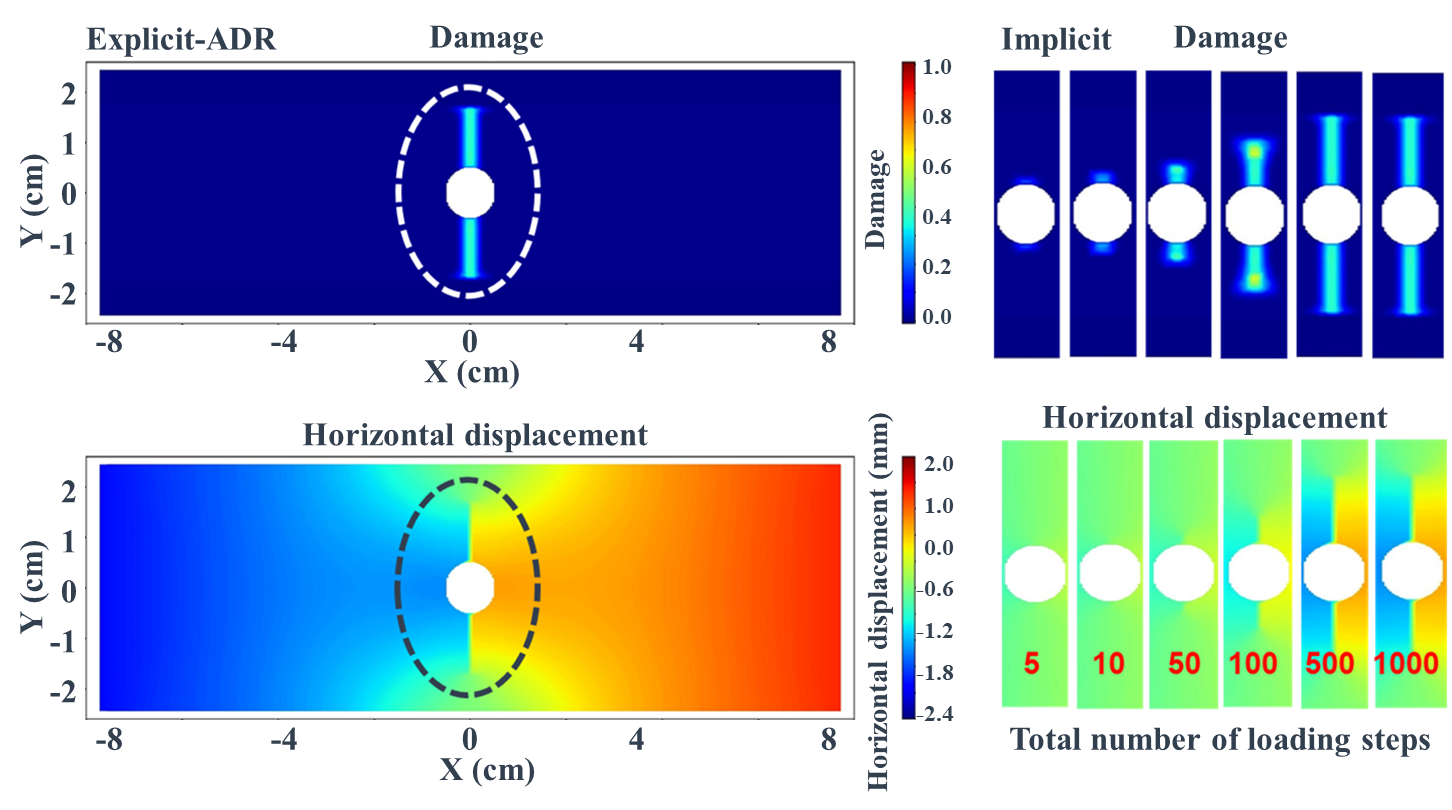}
\caption{The damage and horizontal displacement distribution of a 2D plate under horizontal tension are depicted. The left figure illustrates the results obtained from the explicit-ADR scheme, while the right figure displays the outcomes of the implicit scheme with different total loading steps. It is noteworthy that the right figure only represents the portion enclosed within the dashed boxes of the left figure.}\label{fig6}
\end{figure}

The 2D plate has the initial length $l$ of 0.15m, and the width $w$ of 0.05m. It is uniformly discretized into 7800 PD particles using the regular grid of 150 × 50, with two layers of fictitious thickness 3$\Delta x$ at the displacement boundary on both the left and right sides of the plate. The radius of the central hole is 0.005m. The horizon size is 3.015$\Delta x$. The material is assumed to be isotropic and linear elastic with the Young’s modulus $E$ of 192 GPa and the density of 8000 kg/$m^3$. The boundary and loading conditions are depicted in Figure \ref{fig6}. The left and right ends are subjected to the outward vertical displacement  loadings with a total magnitude of 1.6mm.\par

The convergence criterion remains the same as described in Section \ref{subsec3.1}. Fracture evolution is taken into consideration, and both the implicit and explicit simulations employ the same degradation function to characterize bond damage behavior, as outlined in Eq.(\ref{eq12}), with $s_m = 0.033$, $s_c = 0.066$, and $\beta = 3$.\par
It is worth noting that in this example, the implicit simulation does not converge in one single step but requires multiple loading steps for convergence. For the explicit-ADR, the total displacement loading is attained by applying a small velocity at the boundary over an extended period. The velocity should be small enough to avoid any dynamic effects\cite{RN34}. Kilic and Madenci observed that compared to low loading rates, the results obtained using the explicit-ADR method at high loading rates exhibit significant differences and no longer adhere to the patterns observed at low loading rates. The dynamic effects originate actually from the assumption of small displacements in deriving the ADR formulation. If the loading rate is large, this fundamental assumption of the ADR method is violated, leading to incorrect results. In our case, the velocity magnitude is 1.6×$10^{-7}$ m/iteration and lasts for 10,000 loading steps.\par
\begin{figure}[H]
\centering
\includegraphics[width=0.82\textwidth]{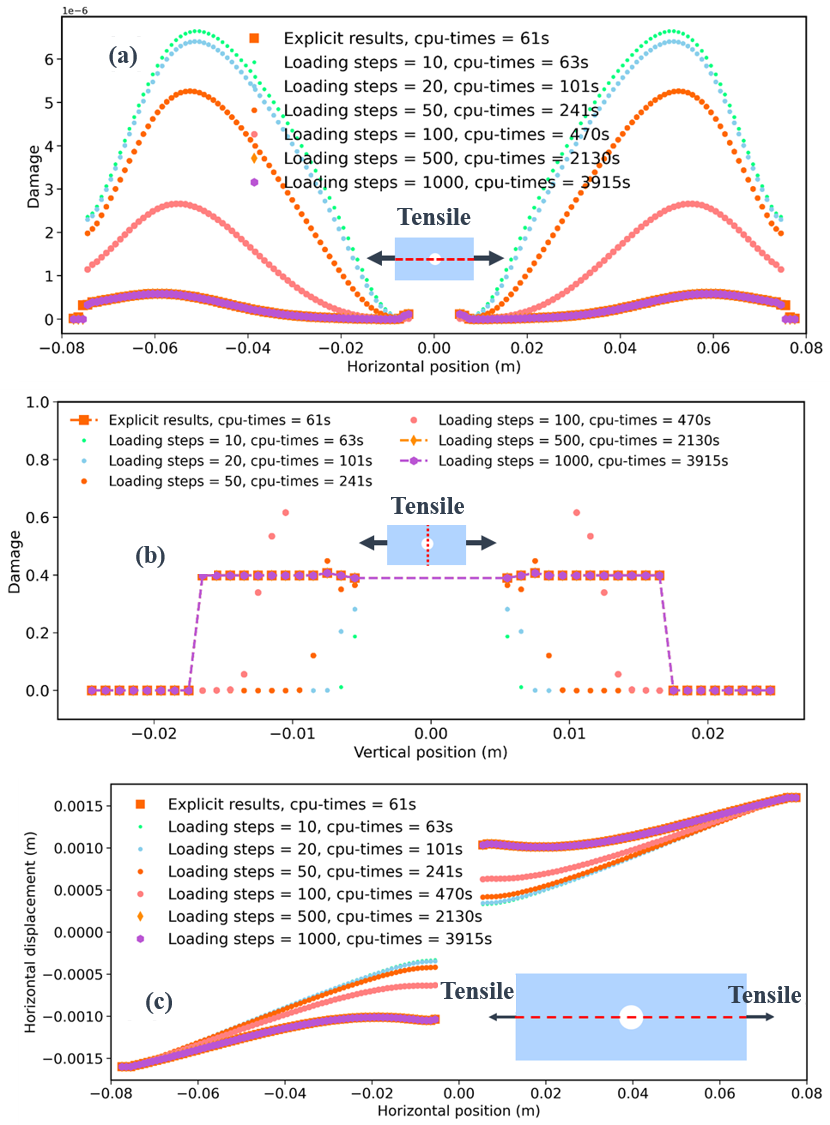}
\caption{This figure illustrates the distribution of damage and horizontal displacement along the specific paths, which are marked as the dashed red lines in each sub-figure. Different colors and markers indicate the change in total loading steps.}\label{fig7}
\end{figure}

Here, we can observe from Figure \ref{fig6} that the results obtained from the implicit schemes vary significantly with the total loading steps. The results between the explicit and implicit simulations are only consistent when the total loading steps is large enough in  the implicit process, e.g., 500 and 1000. To visually observe the differences in the simulation results more intuitively, the damage distribution on the horizontal and vertical central lines are shown in Figure \ref{fig7}.\par

It is evident that both damage and displacement vary significantly with the total loading steps. Moreover, the results obtained from the implicit simulation converge to those from the explicit simulation as the total loading steps increases.\par
As this case involves fracture propagation, the computational costs of the implicit scheme are at least 30 times larger than those of the explicit one. This is attributed to the requirement for a large number of loading steps to ensure computational accuracy.\par
\subsection{Analysis of the implicit scheme  for fracture evolution}
\label{subsec3.3}
Figure \ref{fig7} illustrates that the total loading steps is crucial for computing fracture propagation. In essence, a smaller total number of steps leads to larger loading increments for each step, influencing the bond failure process and consequently leading to different simulations results. Previous studies have suggested that controlling the maximum number of failed bonds at each step is essential for accurate implicit simulations\cite{RN48}. In other words, the smaller the loading increment, the less bond failure occurs at each step, and the more accurate the implicit scheme is for fracture problems. Increasing the total loading steps to reduce the loading increment 
at each step is one of the most effective methods for controlling the maximum number of failed bonds 
at each step, even though it may significantly decrease the efficiency of the implicit scheme. It can be concluded that the implicit scheme is not an efficient choice to simulate fracture problems compared to the explicit-ADR method.\par
When a bond fails, the local forces applied to a PD point experience significant variations in terms of both magnitude and direction, regardless of the type of bond failure,i.e., whether it occurs as a sudden release or through a degradation model.
The differentiation or finite difference in the forces acting on each PD point constitutes the primary element of the Jacobian matrix as depicted in Eq.(\ref{eq6}). If the local forces exhibit pronounced fluctuations, the elements of the Jacobian matrix are inevitably affected 
and undergo significant variations. In the implicit format, the variations of PD forces are reflected in the elements of the Jacobian matrix. \par
\begin{figure}[H]
\centering
\includegraphics[width=0.99\textwidth]{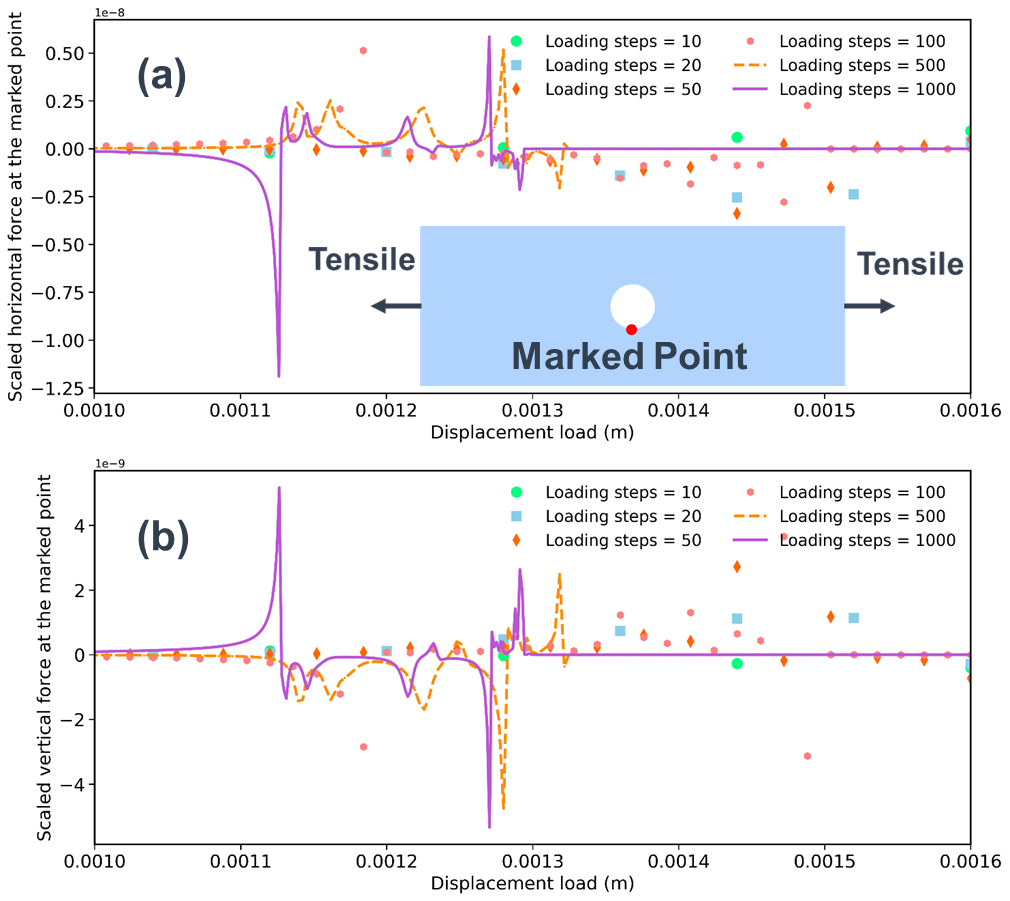}
\caption{Figures (a) and (b) illustrate the variation of  horizontal and vertical forces at the marked points as the total loading steps changes, while the total displacement remains constant. The local force density is obtained during the implicit simulation procedure. Different colors and markers indicate changes in the total loading steps. The scaling factor is denoted as $Q$ ($Q={\Delta x}^3$). The total displacements remain constant at 1.6mm, while the total loading steps vary, ranging from 10 to 2000.}\label{fig8}
\end{figure}

To address the significant variations in local forces during fracturing process, it's necessary to use sufficiently small loading increments, resulting in a large yet essential total loading steps. Figure \ref{fig8} presents the scaled forces of the first-damaged PD point (marked point) during the loading process of the case discussed in Section \ref{subsec3.2}.\par

It is evident that a small total loading steps  smooth out the fluctuations in the local forces and significantly under-predict the variation of local forces once the fractures occur, which fails to describe the interaction forces between the PD points, leading to inaccurate results. Only with a sufficiently large total loading  steps, the variation of local PD forces during fracturing process can be accurately described, making the total loading steps critical for the implicit simulation to calculate fracture evolution. \par
It can be concluded that the implicit scheme is more suitable and effective compared to the explicit-ADR scheme  for quasi-static problems without fractures. This is due to its very fast convergence speed and very few loading steps. However, when fracture evolution is considered, the implicit scheme may require significantly more computational time than the explicit-ADR scheme due to the need for a sufficiently large total loading steps to ensure accuracy. \par
To utilize the efficiency of the implicit scheme for quasi-static problems without fracture and mitigate its disadvantages in simulating fracture
 formation problems, we propose an adaptive strategy for quasi-static fracture formation problems, which will be presented in Section \ref{sec4}.\par
\section{Adaptive strategy for quasi-static fracture problems}
\label{sec4}

\begin{figure}[H]
\centering
\includegraphics[width=0.99\textwidth]{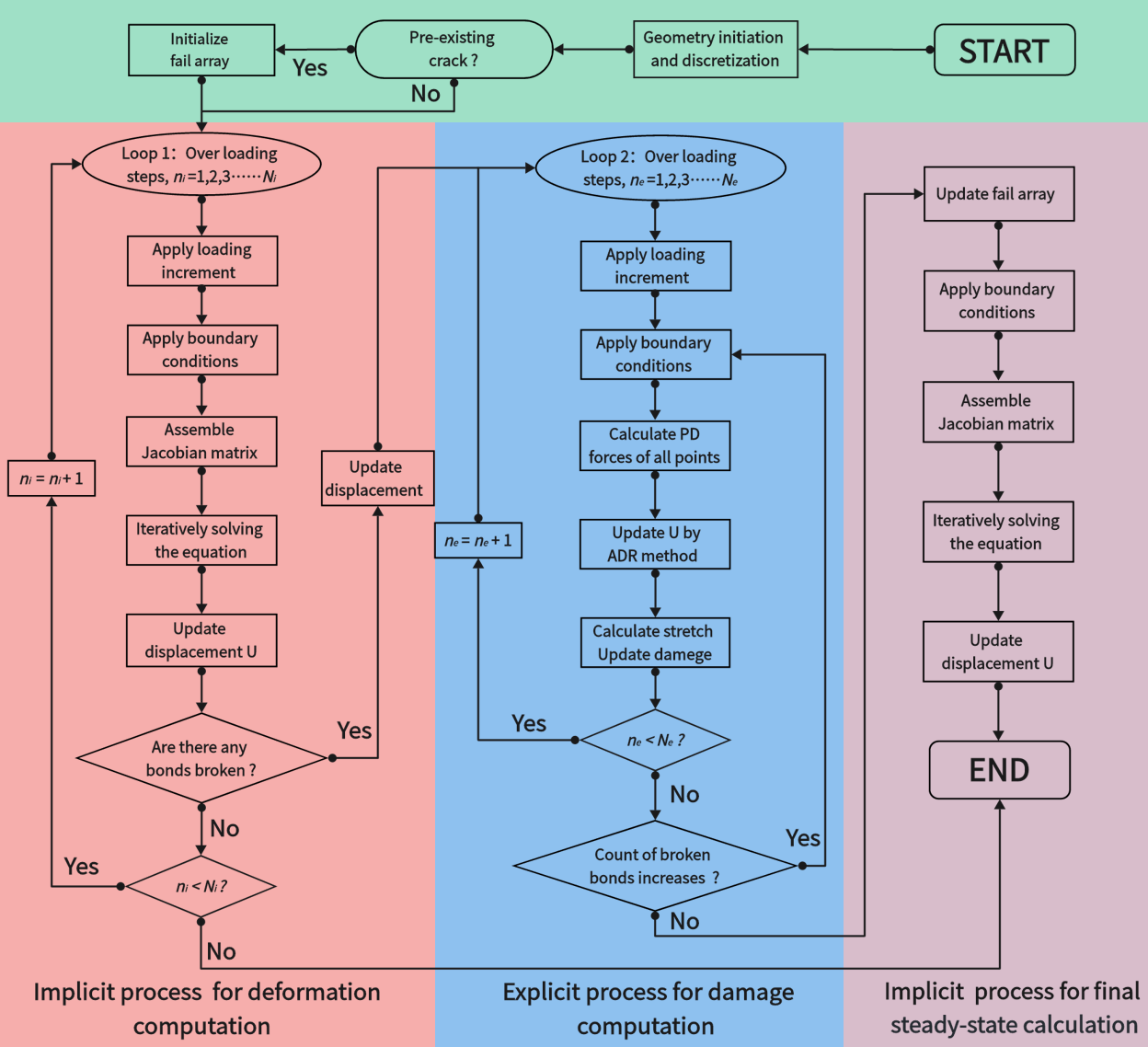}
\caption{The adaptive process flowchart is depicted here. $N_i$ and $N_e$ are the total loading steps of the  implicit and explicit processes, respectively, to reach the total displacement load.}\label{fig9}
\end{figure}
The proposed adaptive scheme is illustrated in Figure \ref{fig9}. The entire process for solving quasi-static fracture problems is divided into three steps. The first step involves the implicit process for deformation computation, extending from the onset of deformation to crack initiation. The second step comprises the explicit-ADR process for damage computation, continuing from the implicit step with inherited displacement and loading information, which proceeds until bond failure ceases. The final step switches back to the implicit process until completion.\par
\begin{figure}[H]
\centering
\includegraphics[width=0.99\textwidth]{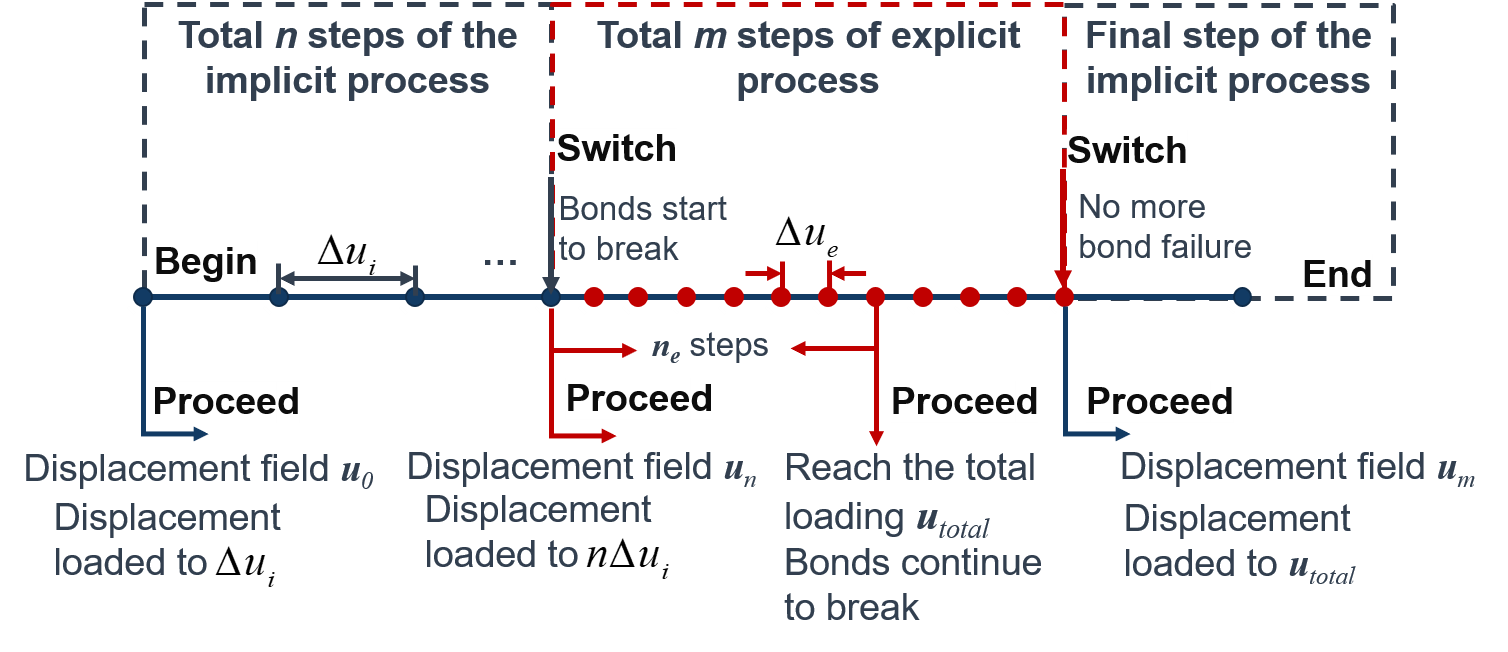}
\caption{The schematic diagram of displacement loading throughout the entire computation process, where  $\Delta u_i, \Delta u_e$ denote the displacement loading increments for the implicit and explicit simulations, respectively, and $u_{total}$ is the total displacement loading. The initial implicit process consists of $n$ steps, while the explicit process consists of $m$ steps, with $N_e$ steps to reach the total displacement load. The relationship between the displacement increments and the total displacement is ${u_{total}} = n \times \Delta u_i + N_e \times \Delta u_e$. }\label{fig10}
\end{figure}

When the tensile stretch of any bond exceeds a certain threshold, denoted as $s_m$, the implicit process terminates and switches to the explicit ADR process. Once the displacement reaches a final state, where no more  bonds fail over an extended period, the explicit ADR process terminates, and the final implicit process resumes. Figure \ref{fig10} illustrates the displacement increments for both the implicit and explicit schemes.\par

It's worth noting that $n$ in Figure \ref{fig10} represents the actual number of computational steps executed in the initial implicit simulation of a general  case, while $N_i$ in Figure \ref{fig9} is the number of steps set before calculation to generate an implicit loading increment $\Delta u_i$ ,  $\Delta u_i = u_{total} / N_i$. In all the simulations, $n$ is not larger than $N_i$.\par 

\section{Numerical validation}
\label{sec5}
We apply the present adaptive and explicit-ADR method to solve four different quasi-static problems. All the  calculations are conducted using Julia programming utilizing a single core of the AMD EPYC 7763 64-Core Processor. The results and the processing times are compared using the same converge criterion.
\subsection{Deformation of undamaged 3D bar under a  transverse loading}
\label{subsec5.1}
\begin{figure}[H]
\centering
\includegraphics[width=0.95\textwidth]{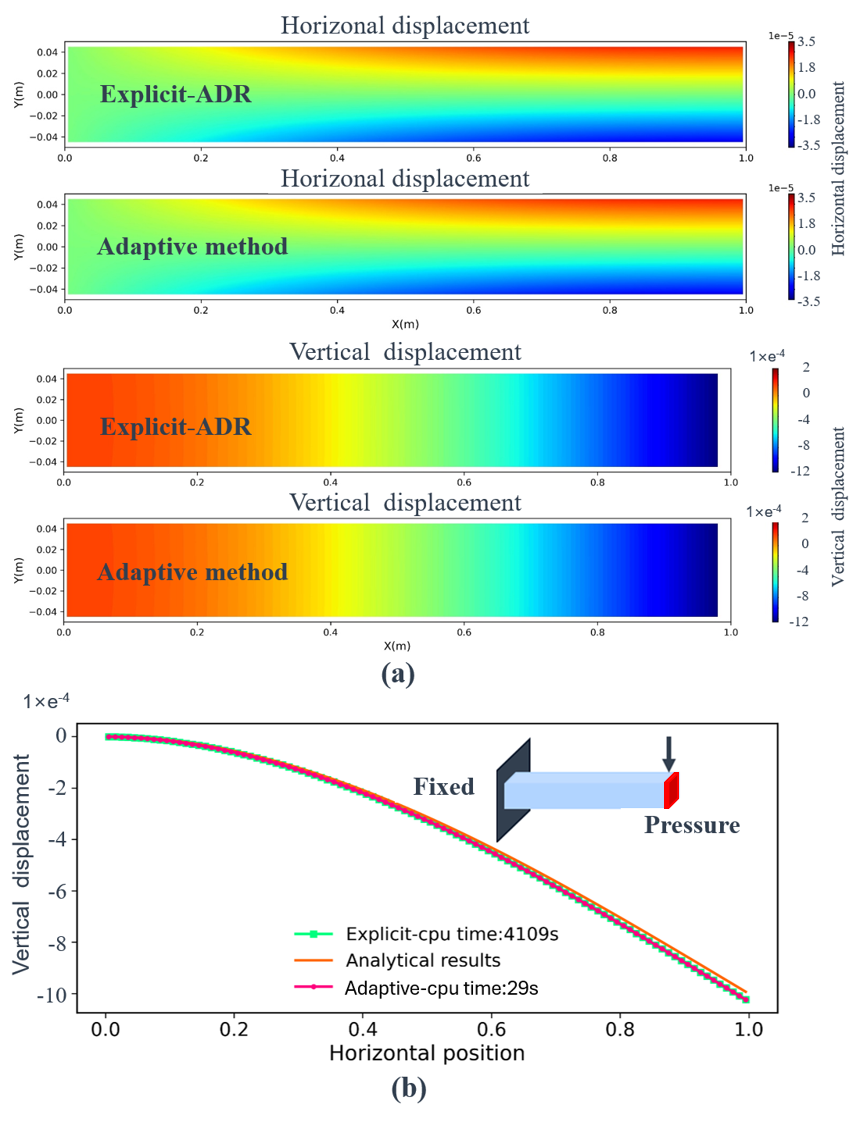}
\caption{Figure (a) shows the contour diagrams of horizontal and vertical displacement obtained by the explicit-ADR and adaptive methods. Figure (b) shows the vertical displacement distribution along the horizontal central line of the bar, and the schematic diagram of geometry, boundary and loading. The loading position is marked in red.}\label{fig11}
\end{figure}

In this benchmark case \cite{RN10}, the 3D bar has the initial length $l$ of 1m, the width $w$ of 1m, and the height $h$ of 0.1m. The domain is uniformly discretized into 10,300 particles using the regular grid of  100 × 10 × 10, with an additional layer of fictitious thickness 3$\Delta x$ at the fixed boundary on the left side of the bar. The horizon size is 3.015$\Delta x$. The material is assumed to be isotropic and linear elastic with the Young’s modulus $E$ of 200 GPa and the density $\rho$ of 7850 kg/$m^3$. The boundary and loading conditions are depicted in Figure \ref{fig11}(b). The left end of the bar is fixed, while the right end is subjected to a downward vertical loading of 5000 N. The loading boundary is located on the outermost layer of the right end of the bar, with a boundary layer volume of 10×10×$(\Delta x)^3$. The loading is applied from the beginning as a body force density vector to the PD points within the loading boundary layer, with a value of 5×$10^7$ N/$m^3$. \par

In this case, both the adaptive and explicit-ADR methods directly utilize the total pressure loading as the boundary condition and adhere to the same convergence criterion, as described in Eq.(\ref{eq16}), where $e_0$ is set to be  1×$10^{-9}$. The adaptive method achieves convergence in just one step, while the explicit one requires approximately 20,000 steps. \par
We found that the results of the present adaptive method are consistent with those of the explicit scheme reported in the literature\cite{RN10}, with an acceleration ratio of approximately 141.7.\par

\subsection{Damage evolution of a 2D square plate with a central hole under tensile loadings}
\label{subsec5.2}
This case is similar to the benchmark example of the citation \cite{RN46}, the 2D square plate has the initial length $l$ of 0.05m. It is uniformly discretized into 2900 particles using the regular grid of 50 × 50, with two additional layers of fictitious thickness 3$\Delta x$ at the displacement boundaries on both the top and bottom sides of the plate. The radius of the central hole is 0.005m. The horizon size is 3.015$\Delta x$. The material is assumed to be isotropic and linear elastic with the Young’s modulus $E$ of 192 GPa and the density $\rho$ of 8000 kg/$m^3$. The bond failure criterion is based on the degradation form as outlined in Eq.(\ref{eq12}), with $s_m$  = 0.015, $s_c$ = 0.02, and $\beta$ = 3. The boundary and loading conditions are depicted in Figure \ref{fig12}(b). The upper and lower ends of the plate are subjected to the outward vertical displacement loadings with a magnitude of 0.275 mm.\par

\begin{figure}[H]
\centering
\includegraphics[width=0.9\textwidth]{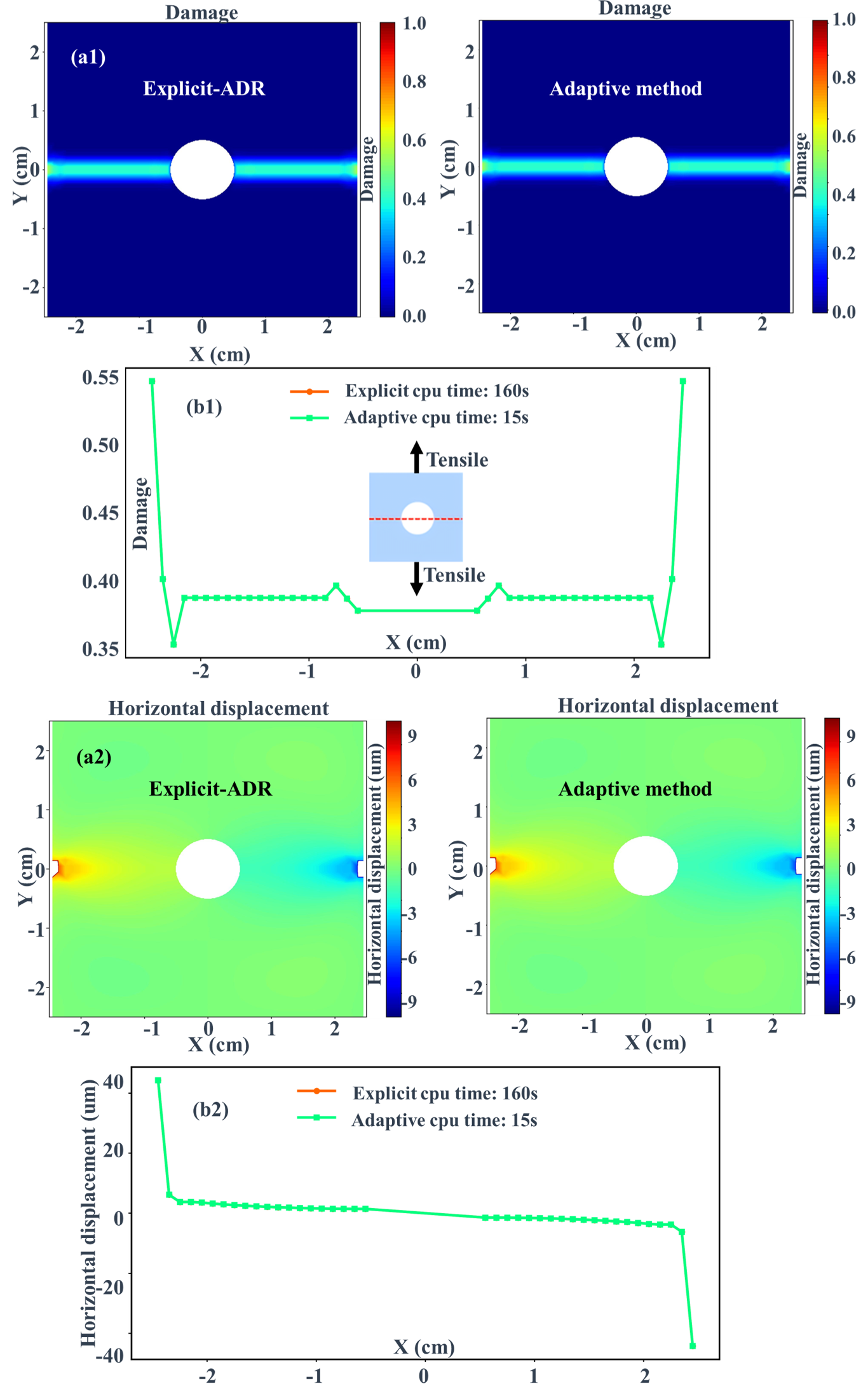}
\label{fig12a} 
\end{figure}

\begin{figure}[H]
\centering
\includegraphics[width=0.99\textwidth]{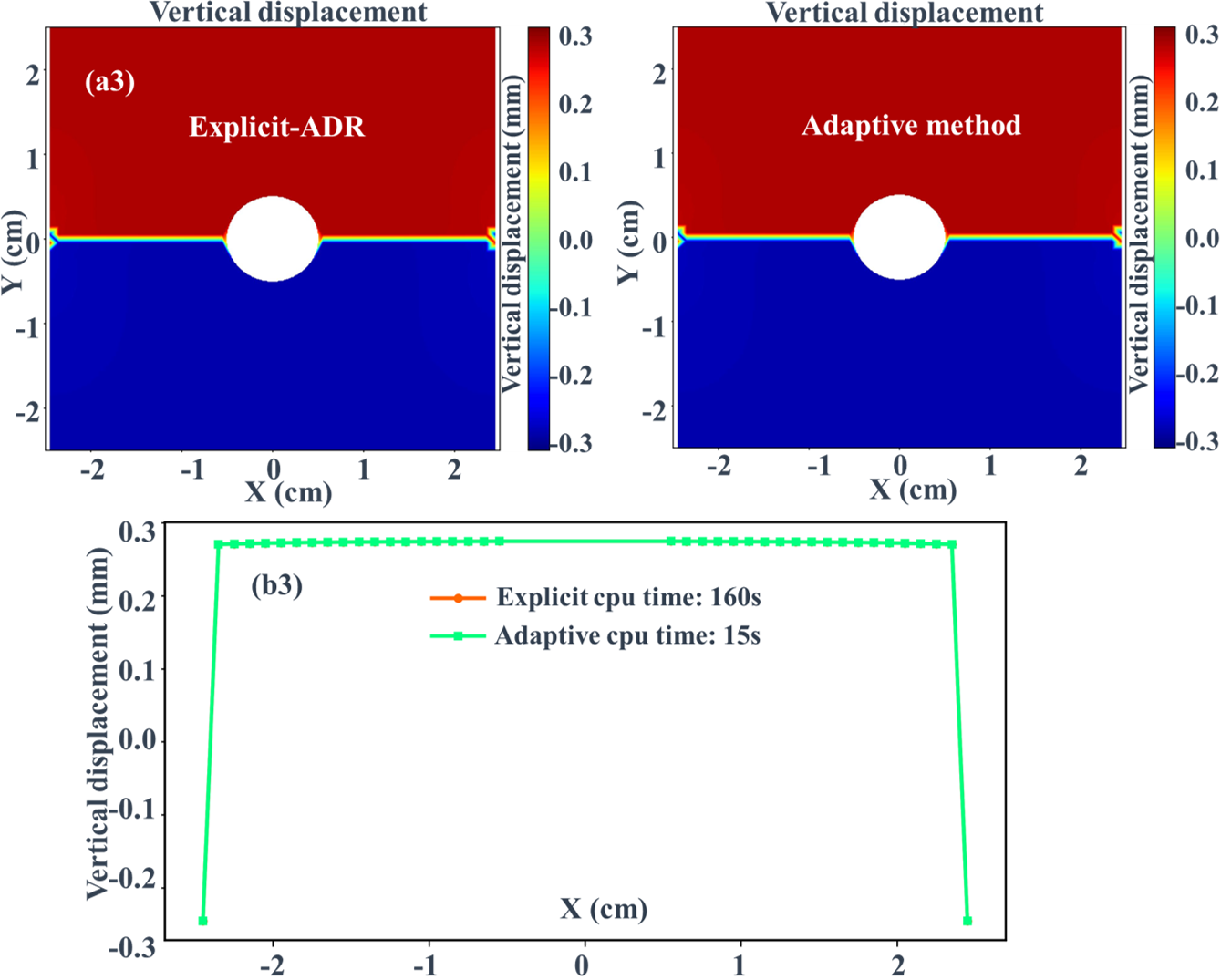}
\caption{The schematic diagram of geometry, boundary and loading is shown in Figure (b1). Figures (a1-a3) display the contour diagrams illustrating the distribution of damage, horizontal displacement, and vertical displacement. Figures (b1-b3) show the distribution of damage, horizontal displacement, and vertical displacement along the horizontal central line of the bar, indicated by the red dash line in Figure (a1).}
\label{fig12}
\end{figure}

In this case, the total displacement loading is applied through both implicit and explicit loading processes in the adaptive method. The maximum number of steps for the implicit and explicit loading processes is set to be 3 and 180, respectively. This implies that the loading rate for deformation computation is \(2.75 \times 10^{-4}/3\) m/iteration, while the loading rate for damage evolution computation is \((1 - n/3) \times 2.75 \times 10^{-4}/180\) m/iteration, where \(n\) is the step at which damage initiation occurs, as depicted in Figure \ref{fig10}. 
As for the explicit-ADR scheme, the loading rate is set to be 2.75×$10^{-7}$ m/iteration. In addition, both the adaptive and explicit methods use the same convergence criterion, as described in Eq.(\ref{eq16}), where $e_0$ is set to be  1×$10^{-9}$.\par
The adaptive computation results are in excellent agreement with those obtained by the explicit-ADR simulation, with an acceleration ratio of 10.7.\par

\subsection{ Damage evolution of a 2D plate subjected to  three point bending}
\label{subsec5.3}
\begin{figure}[H]
\centering
\includegraphics[width=0.95\textwidth]{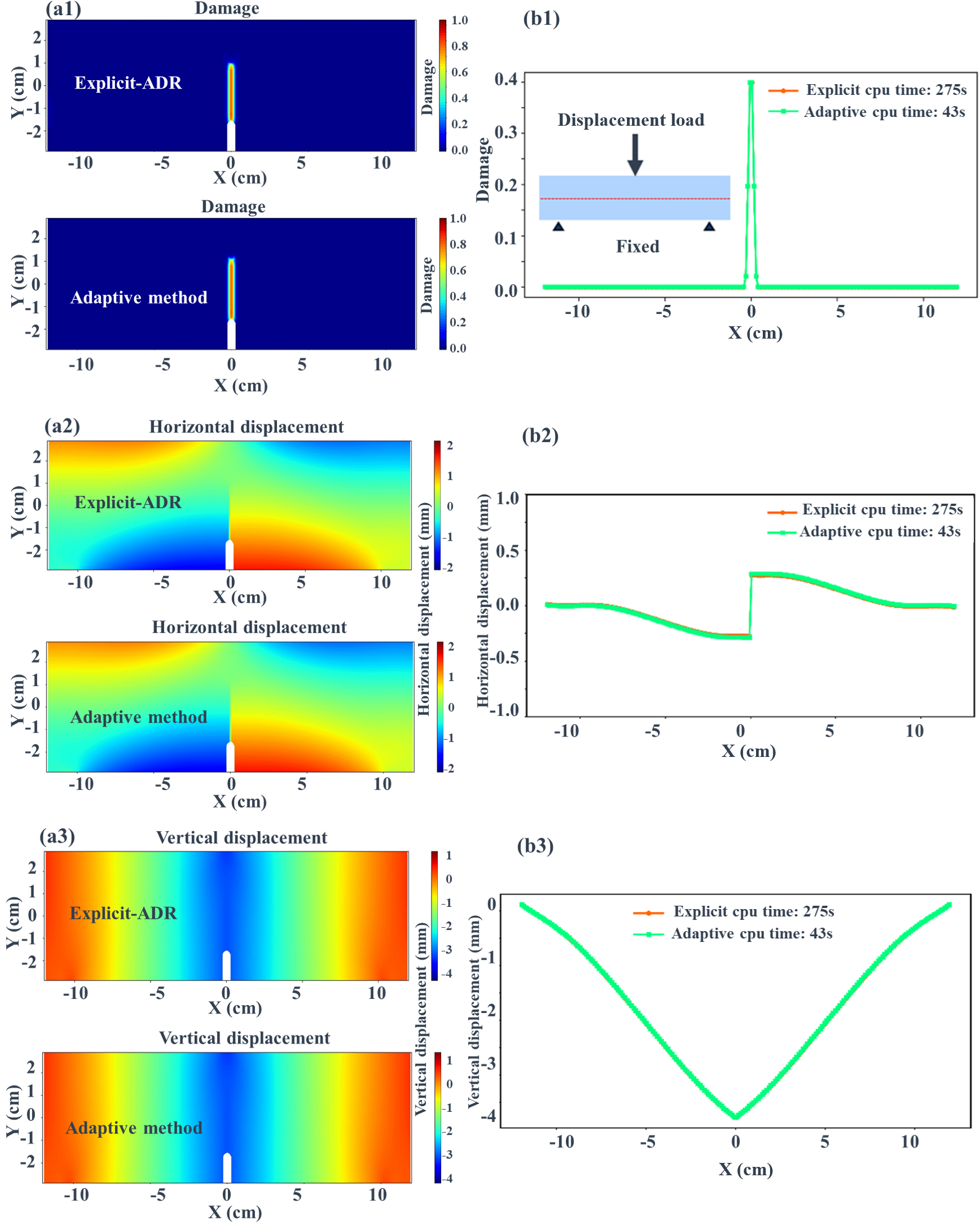}
\caption{The schematic diagram of geometry, boundary and loading is shown in Figure (b1). Figures (a1-a3) display the contour diagrams for the distribution of damage, horizontal displacement, and vertical displacement, respectively. Figures (b1-b3) show the damage, horizontal displacement, and vertical displacement along the horizontal central line of the plate, respectively.}\label{fig13}
\end{figure}

In this case, the 2D rectangular plate has the initial length $l$ of 0.24 m, and the width $w$ of 0.06 m. It is uniformly discretized into 10,048 PD particles using the regular grid of 200× 50, with three additional areas of fictitious thickness 3$\Delta x$, and the span 4$\Delta x$ at the displacement boundary on the three point for loading. The pre-existing notch is located at the center of the lower end, oriented vertically, with the length of 0.3$w$ and the span of 2$\Delta x$. The horizon size is 3.015$\Delta x$. The material is assumed to be isotropic and linear elastic with the Young’s modulus $E$ of 200 GPa and the density $\rho$ of 8000 kg/$m^3$. The bond failure criterion is based on the degradation form as described in Eq.(\ref{eq12}), with $s_m$ = 0.016, $s_c$ = 0.02, and $\beta$ = 3. The boundary and loading conditions are depicted in Figure \ref{fig13} (b1). \par

As for the three point bending test, two points at the lower end are constrained from the vertical movement, each with a horizontal distance of 0.105 m from the center. Additionally, one last point located at the center of the upper end is subjected to the downward vertical displacement loading with a magnitude of 4 mm.\par
Here, the maximum number of steps for the implicit and explicit loading processes is set to be 5 and 1000 for the adaptive method, respectively. This implies that the loading rate for deformation computation is 4×$10^{-3}$/5 m/iteration, and for damage evolution computation is (1 – $n$/5)×4×$10^{-3}$/1000 m/iteration, where $n$ is the step at which damage initiation occurs, as depicted in Figure \ref{fig10}. As for the explicit-ADR scheme, the loading rate is set to be  1×$10^{-7}$ m/iteration. In addition, both the adaptive and explicit methods use the same convergence criterion, as described in Eq.(\ref{eq16}), where $e_0$ is set to be 1×$10^{-9}$.\par
The results of the current adaptive method are in excellent agreement with those obtained by the explicit scheme, exhibiting an acceleration ratio of around 6.4.

\subsection{Damage evolution of a 2D plate with multiple holes under a concentrated tensile loading}
\label{subsec5.4}
\begin{figure}[H]
\centering
\includegraphics[width=0.35\textwidth]{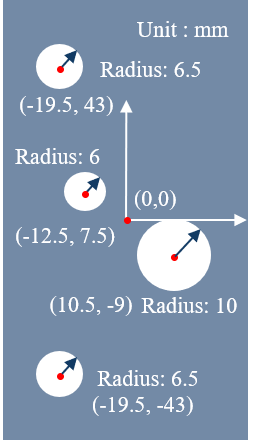}
\caption{The detailed geometric layout, where the values in parentheses represent the coordinates of the hole centers.}\label{fig14}
\end{figure}

In this case, the 2D plate has the initial length $l$ of 0.065 m, the width $w$ of 0.12 m, as described in \cite{RN54}. It is uniformly discretized into about 7800 particles using the regular grid of 65 × 120. There are four holes in the plate. The radius of the upper left and lower left holes for loading is 6.5 mm, while the radius of the smaller hole in the middle is 6 mm, and the radius of the larger one is 10 mm. The geometric layout is shown in Figure \ref{fig14}.\par

Note, the horizon size is set to be  8.015$\Delta x$ for fracture evolution calculations to eliminate dependence on the regular distribution of PD points during crack propagation\cite{RN55}. For the adaptive method, a horizon size of 3.015$\Delta x$ is used throughout the material during the implicit simulation of the deformation stage.\par
\begin{figure}[H]
\centering
\includegraphics[width=0.91\textwidth]{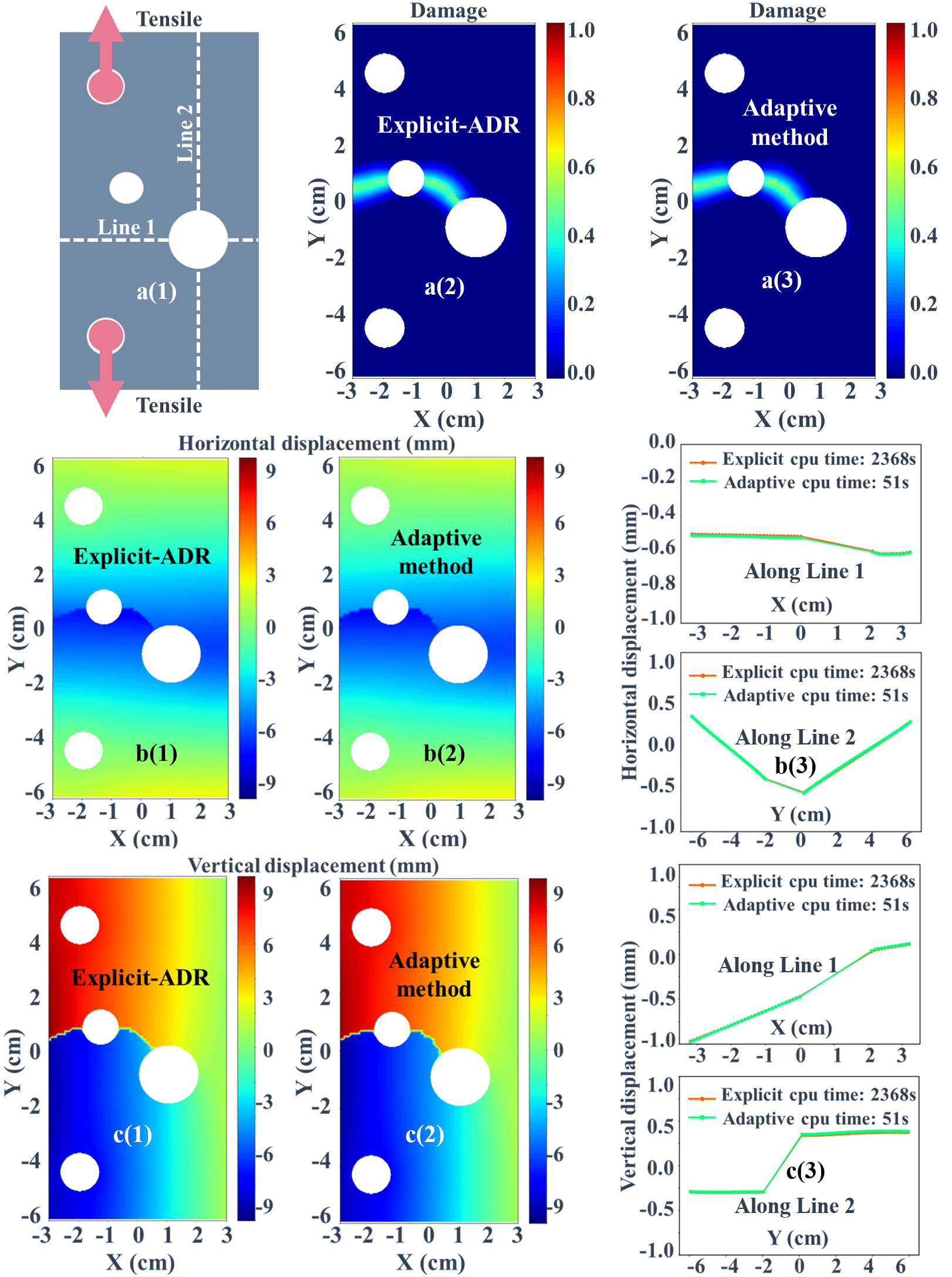}
\caption{Comparison of the explicit and adaptive methods. Figures (a2-a3) show the damage distribution, Figures (b1-b2) display the horizontal displacement, and Figures (c1-c2) present the vertical displacement. Figures (b3) and (c3) illustrate the horizontal and vertical displacement along Line 1 and Line 2, which are marked in Figure (a1) by the white dashed lines.}\label{fig15}
\end{figure}

In all the simulations, the linear isotropic elastic model is applied, with the Young’s modulus $E$ of 200 GPa and the density $\rho$ of 8000 kg/$m^3$. The bond failure criterion follows Eq.(\ref{eq12}), with $s_m$ = 0.016, $s_c$ = 0.02, and $\beta$ = 3. The boundary and loading conditions are shown in Figure \ref{fig15}(a1), where the holes on the left side experience the outward vertical displacement loadings of 8.0×$10^{-4}$m.

Here, the maximum number of steps for the implicit and explicit loading processes is set to be 5 and 1500 for the adaptive method, respectively. This implies that the loading rate for deformation computation is 8.0×$10^{-4}$/5 m/iteration, and for damage evolution computation is (1 –$n$/5)×8.0×$10^{-4}$/1500 m/iteration,  respectively, where $n$ is the step at which damage initiation occurs, as depicted in Figure \ref{fig10}. As for the explicit-ADR scheme, the loading rate is set to be 1×$10^{-8}$ m/iteration. Both the adaptive and explicit methods use the same convergence criterion with $e_0$ = 1×$10^{-9}$ as described in Eq.(\ref{eq16}).\par
The results of the current adaptive method are in excellent agreement with those obtained by the explicit scheme, exhibiting an acceleration ratio of around 46.4.\par
In summary, the table below lists the computational time consumed for the above four test cases. All the calculations were conducted, using Julia programming, on a single core of the AMD EPYC 7763 64-Core Processor.\par
\begin{figure}[H]
\centering
\includegraphics[width=0.75\textwidth]{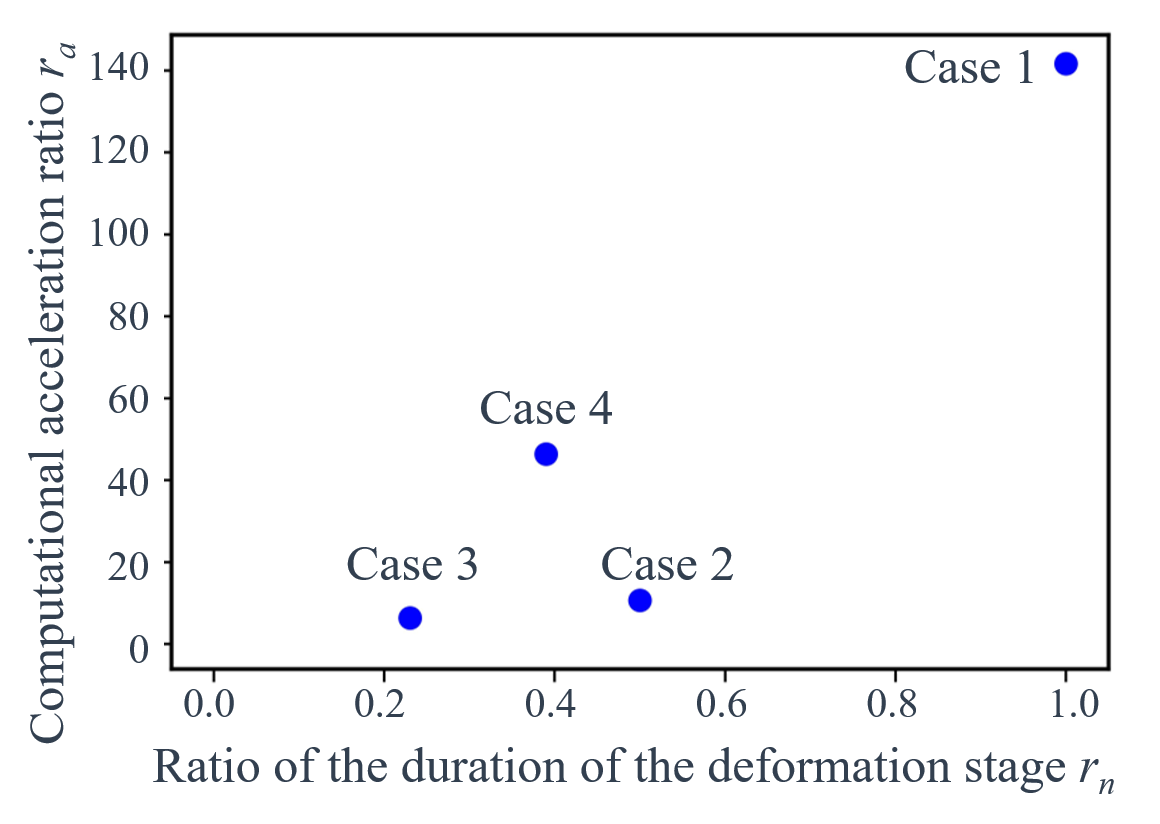}
\caption{ The relationship between $r_n$ and $r_a$ for the four test cases.}\label{fig16}
\end{figure}

\setcounter{table}{0}
\begin{table}[h]
    \centering
    \caption{Summary of computational time}
    \label{table:cases}
    \begin{tabular}{cccc}
        \toprule
        Case Number & \begin{tabular}[c]{@{}c@{}}explicit-ADR\\ method (s)\end{tabular} & \begin{tabular}[c]{@{}c@{}}Adaptive\\ method (s)\end{tabular} & \begin{tabular}[c]{@{}c@{}}Acceleration ratio of\\ adaptive to explicit-ADR method\end{tabular} \\
        \midrule
        1 & 4109 & 29 & 141.7 \\
        2 & 160 & 15 & 10.7 \\
        3 & 275 & 43 & 6.4 \\
        4 & 2368 & 51 & 46.4 \\
        \bottomrule
    \end{tabular}
\end{table}
It can be observed from the above simulations that the present adaptive method can efficiently solve quasi-static fracture problems, akin to the explicit-ADR scheme, which reduces computational costs significantly compared to the explicit scheme. We observe that the acceleration ratio $r_a$ of the adaptive method in comparison with the explicit-ADR scheme for the above four cases ranging from 6.4 to 141.7. In the entire simulation process of 
the 4 cases, including deformation and fracture propagation, we calculate the ratio $r_n$ of the number of steps from the beginning to fracture initiation to the total number of steps required for convergence to the  specified precision, and reveal the relationship between acceleration ratios $r_a$ and $r_n$ as shown in Figure  \ref{fig16}. \par

It is evident that $r_n$ plays a significant role in the acceleration performance of the adaptive method: the larger $r_n$, the greater proportion of the deformation process in the entire process, and the better acceleration ratio due to more significant contribution of the implicit scheme. The result of Case 4 is unexpected, as its acceleration ratio $r_a$ is much larger than anticipated. The reason for this discrepancy is that, in Case 4, the acceleration effect stems not only from the greater proportion of the deformation process but also from the much smaller horizon size of 3.015$\Delta x$ for the deformation stage of the adaptive method compared to 8.015$\Delta x$ for the explicit scheme. Furthermore, due to the larger horizon size, which would result in a significant computational burden, the acceleration effect of the latter is more pronounced.\par

\section{Concluding Remarks}
\label{sec6}
   In this study, we have developed an efficient implicit BBPD model for solving quasi-static problems based on the full nonlinear equilibrium equation with a degradation bond failure function. We compared the computational efficiency of the explicit-ADR method with the newly-developed implicit method and further proposed an adaptive explicit-implicit method to enhance computational efficiency for solving quasi-static problems with PD. Three main conclusions drawn from this work:\par
(1)	The Jacobian matrix based on the full nonlinear BBPD, derived from the NR procedure, exhibits typical properties such as sparsity, symmetry, and  non-singularity, indicating good solvability. \par
(2)	Regarding the computational efficiency of solving quasi-static problems using PD, the present implicit method performs well in computing damage-free deformations but struggles with damage evolution. This is because a sufficiently large number of loading steps are required for the implicit simulations to accurately compute damage evolution.\par
(3) The proposed adaptive method utilizes the advantages of the implicit and explicit methods to improve overall computational efficiency. The accuracy and efficiency of our adaptive method are examined in the four test cases, revealing that the greater the proportion of the deformation process in the entire process, the more pronounced computational  acceleration of the adaptive method. \par
Meanwhile, three key points are worthwhile for in-depth research in the future.\par
  (1) No bonds fail during deformation, but the implicit simulation needs to calculate bond stretch for all the PD points at each step to decide whether to switch to the explicit scheme, which leads to substantial computational burden. The switching criteria may therefore be refined to avoid calculating bond stretch during the deformation stage, thus further reducing computational costs.\par
  (2) The conjugate gradient method is used to solve the equations. Multiple numerical techniques, including preconditioning, may be exploited to expedite computations for large linear systems and further optimize the adaptive method.\par
  (3) While the present explicit-implicit adaptive scheme is demonstrated for BBPD, it can be straightforwardly  extended to SBPD.

\section*{Acknowledgment}
The authors acknowledge the financial support of the National Science Foundation of China (No.12302381) and the computing resources provided by Hefei Advanced Computing Center.

\bibliographystyle{elsarticle-num}
\bibliography{hybrid} 

\begin{thebibliography}{10}
\expandafter\ifx\csname url\endcsname\relax
  \def\url#1{\texttt{#1}}\fi
\expandafter\ifx\csname urlprefix\endcsname\relax\def\urlprefix{URL }\fi
\expandafter\ifx\csname href\endcsname\relax
  \def\href#1#2{#2} \def\path#1{#1}\fi

\bibitem{RN1}
S.~A. Silling, Reformulation of elasticity theory for discontinuities and
  long-range forces, Journal of the Mechanics and Physics of Solids 48~(1)
  (2000) 175--209.
\newblock \href {https://doi.org/10.1016/S0022-5096(99)00029-0}
  {\path{doi:10.1016/S0022-5096(99)00029-0}}.

\bibitem{RN2}
S.~A. Silling, E.~Askari, Peridynamic modeling of impact damage, in: ASME
  pressure vessels and piping conference, Vol. 46849, 2004, pp. 197--205.

\bibitem{RN3}
F.~Bobaru, S.~A. Silling, H.~Jiang, Peridynamic fracture and damage modeling of
  membranes and nanofiber networks., Tech. rep., Sandia National Lab.(SNL-NM),
  Albuquerque, NM (United States) (2004).

\bibitem{RN4}
S.~A. Silling, E.~Askari, A meshfree method based on the peridynamic model of
  solid mechanics, Computers \& structures 83~(17-18) (2005) 1526--1535.
\newblock \href {https://doi.org/10.1016/j.compstruc.2004.11.026}
  {\path{doi:10.1016/j.compstruc.2004.11.026}}.

\bibitem{RN5}
P.~N. Demmie, S.~A. Silling, An approach to modeling extreme loading of
  structures using peridynamics, Journal of Mechanics of Materials and
  Structures 2~(10) (2007) 1921--1945.
\newblock \href {https://doi.org/10.2140/jomms.2007.2.1921}
  {\path{doi:10.2140/jomms.2007.2.1921}}.

\bibitem{RN6}
W.~Gerstle, N.~Sau, S.~Silling, Peridynamic modeling of concrete structures,
  Nuclear Engineering and Design 237~(12-13) (2007) 1250--1258.
\newblock \href {https://doi.org/10.1016/j.nucengdes.2006.10.002}
  {\path{doi:10.1016/j.nucengdes.2006.10.002}}.

\bibitem{RN7}
S.~A. Silling, M.~Epton, O.~Weckner, J.~Xu, E.~Askari, Peridynamic states and
  constitutive modeling, Journal of Elasticity 88~(2) (2007) 151--184.
\newblock \href {https://doi.org/10.1007/s10659-007-9125-1}
  {\path{doi:10.1007/s10659-007-9125-1}}.

\bibitem{RN9}
S.~Silling, R.~Lehoucq, Peridynamic theory of solid mechanics, in: H.~Aref,
  E.~van~der Giessen (Eds.), Advances in Applied Mechanics, Vol.~44 of Advances
  in Applied Mechanics, Elsevier, 2010, pp. 73--168.
\newblock \href {https://doi.org/https://doi.org/10.1016/S0065-2156(10)44002-8}
  {\path{doi:https://doi.org/10.1016/S0065-2156(10)44002-8}}.

\bibitem{RN10}
E.~Madenci, E.~Oterkus, Peridynamic Theory and Its Applications, Springer New
  York, 2014.

\bibitem{RN11}
Q.~Tong, S.~Li, Multiscale coupling of molecular dynamics and peridynamics,
  Journal of the Mechanics and Physics of Solids 95 (2016) 169--187.
\newblock \href {https://doi.org/10.1016/j.jmps.2016.05.032}
  {\path{doi:10.1016/j.jmps.2016.05.032}}.

\bibitem{RN12}
A.~Javili, R.~Morasata, E.~Oterkus, S.~Oterkus, Peridynamics review,
  Mathematics and Mechanics of Solids 24~(11) (2018) 3714--3739.
\newblock \href {https://doi.org/10.1177/1081286518803411}
  {\path{doi:10.1177/1081286518803411}}.

\bibitem{RN13}
A.~Pirzadeh, F.~Dalla~Barba, F.~Bobaru, L.~Sanavia, M.~Zaccariotto,
  U.~Galvanetto, Elastoplastic peridynamic formulation for materials with
  isotropic and kinematic hardening, Engineering with Computers 40~(4) (2024)
  2063--2082.
\newblock \href {https://doi.org/10.1007/s00366-024-01943-x}
  {\path{doi:10.1007/s00366-024-01943-x}}.

\bibitem{RN14}
F.~Han, G.~Lubineau, Y.~Azdoud, Adaptive coupling between damage mechanics and
  peridynamics: A route for objective simulation of material degradation up to
  complete failure, Journal of the Mechanics and Physics of Solids 94 (2016)
  453--472.
\newblock \href {https://doi.org/10.1016/j.jmps.2016.05.017}
  {\path{doi:10.1016/j.jmps.2016.05.017}}.

\bibitem{RN15}
A.~R. Torabi, H.~R. Majidi, H.~Amani, J.~Akbardoost, Implementation of xfem for
  fracture prediction of vo-notched brittle specimens, European Journal of
  Mechanics - A/Solids 81 (2020) 103970.
\newblock \href {https://doi.org/10.1016/j.euromechsol.2020.103970}
  {\path{doi:10.1016/j.euromechsol.2020.103970}}.

\bibitem{RN16}
X.~Zhou, J.~Chen, Extended finite element simulation of step-path brittle
  failure in rock slopes with non-persistent en-echelon joints, Engineering
  Geology 250 (2019) 65--88.
\newblock \href {https://doi.org/10.1016/j.enggeo.2019.01.012}
  {\path{doi:10.1016/j.enggeo.2019.01.012}}.

\bibitem{RN17}
N.~Moës, J.~Dolbow, T.~Belytschko, A finite element method for crack growth
  without remeshing, International Journal for Numerical Methods in Engineering
  46~(1) (1999) 131--150.
\newblock \href
  {https://doi.org/10.1002/(SICI)1097-0207(19990910)46:1<131::AID-NME726>3.0.CO;2-J}
  {\path{doi:10.1002/(SICI)1097-0207(19990910)46:1<131::AID-NME726>3.0.CO;2-J}}.

\bibitem{RNzy1}
Y.~Zhang, S.~Haeri, G.~Pan, Y.~Zhang, Strongly coupled peridynamic and lattice
  boltzmann models using immersed boundary method for flow-induced structural
  deformation and fracture, Journal of Computational Physics 435 (2021) 110267.
\newblock \href {https://doi.org/10.1016/j.jcp.2021.110267}
  {\path{doi:10.1016/j.jcp.2021.110267}}.

\bibitem{RNzy2}
Y.~Zhang, G.~Pan, Y.~Zhang, S.~Haeri, A multi-physics peridynamics-dem-ib-clbm
  framework for the prediction of erosive impact of solid particles in viscous
  fluids, Computer Methods in Applied Mechanics and Engineering 352 (2019)
  675--690.
\newblock \href {https://doi.org/10.1016/j.cma.2019.04.043}
  {\path{doi:10.1016/j.cma.2019.04.043}}.

\bibitem{RN18}
W.~Hu, Y.~D. Ha, F.~Bobaru, S.~A. Silling, The formulation and computation of
  the nonlocal j-integral in bond-based peridynamics, International Journal of
  Fracture 176~(2) (2012) 195--206.
\newblock \href {https://doi.org/10.1007/s10704-012-9745-8}
  {\path{doi:10.1007/s10704-012-9745-8}}.

\bibitem{RN19}
S.~A. Silling, M.~D’Elia, Y.~Yu, H.~You, M.~Fermen-Coker, Peridynamic model
  for single-layer graphene obtained from coarse-grained bond forces, Journal
  of Peridynamics and Nonlocal Modeling 5~(2) (2023) 183--204.
\newblock \href {https://doi.org/10.1007/s42102-021-00075-w}
  {\path{doi:10.1007/s42102-021-00075-w}}.

\bibitem{RN20}
Y.~Wang, X.~Zhou, M.~Kou, A coupled thermo-mechanical bond-based peridynamics
  for simulating thermal cracking in rocks, International Journal of Fracture
  211~(1-2) (2018) 13--42.
\newblock \href {https://doi.org/10.1007/s10704-018-0273-z}
  {\path{doi:10.1007/s10704-018-0273-z}}.

\bibitem{RN21}
Y.~Wang, X.~Zhou, Y.~Wang, Y.~Shou, A 3-d conjugated bond-pair-based
  peridynamic formulation for initiation and propagation of cracks in brittle
  solids, International Journal of Solids and Structures 134 (2018) 89--115.
\newblock \href {https://doi.org/10.1016/j.ijsolstr.2017.10.022}
  {\path{doi:10.1016/j.ijsolstr.2017.10.022}}.

\bibitem{RN22}
E.~Madenci, M.~Dorduncu, N.~Phan, X.~Gu, Weak form of bond-associated
  non-ordinary state-based peridynamics free of zero energy modes with uniform
  or non-uniform discretization, Engineering Fracture Mechanics 218 (2019)
  106613.
\newblock \href {https://doi.org/10.1016/j.engfracmech.2019.106613}
  {\path{doi:10.1016/j.engfracmech.2019.106613}}.

\bibitem{RN23}
X.-P. Zhou, X.-L. Yu, A vector form conjugated-shear bond-based peridynamic
  model for crack initiation and propagation in linear elastic solids,
  Engineering Fracture Mechanics 256 (2021) 107944.
\newblock \href {https://doi.org/10.1016/j.engfracmech.2021.107944}
  {\path{doi:10.1016/j.engfracmech.2021.107944}}.

\bibitem{RN24}
P.~Nikolaev, M.~Sedighi, A.~P. Jivkov, L.~Margetts, Analysis of heat transfer
  and water flow with phase change in saturated porous media by bond-based
  peridynamics, International Journal of Heat and Mass Transfer 185 (2022)
  122327.
\newblock \href {https://doi.org/10.1016/j.ijheatmasstransfer.2021.122327}
  {\path{doi:10.1016/j.ijheatmasstransfer.2021.122327}}.

\bibitem{RN25}
S.~Li, H.~Lu, X.~Huang, R.~Qin, J.~Mao, Sensitivity analysis of notch shape on
  brittle failure by using uni-bond dual-parameter peridynamics, Engineering
  Fracture Mechanics 291 (2023) 109566.
\newblock \href {https://doi.org/10.1016/j.engfracmech.2023.109566}
  {\path{doi:10.1016/j.engfracmech.2023.109566}}.

\bibitem{RN26}
A.~Masoumi, M.~Salehi, M.~Ravandi, A modified bond-based peridynamic model
  without limitations on elastic properties, Engineering Analysis with Boundary
  Elements 149 (2023) 261--281.
\newblock \href {https://doi.org/10.1016/j.enganabound.2023.01.030}
  {\path{doi:10.1016/j.enganabound.2023.01.030}}.

\bibitem{RN27}
G.~Ongaro, A.~Shojaei, F.~Mossaiby, A.~Hermann, C.~J. Cyron, P.~Trovalusci,
  Multi-adaptive spatial discretization of bond-based peridynamics,
  International Journal of Fracture 244~(1) (2023) 1--24.
\newblock \href {https://doi.org/10.1007/s10704-023-00709-8}
  {\path{doi:10.1007/s10704-023-00709-8}}.

\bibitem{RN28}
H.~Ren, X.~Zhuang, X.~Fu, Z.~Li, T.~Rabczuk, Bond-based nonlocal models by
  nonlocal operator method in symmetric support domain, Computer Methods in
  Applied Mechanics and Engineering 418 (2024) 116230.
\newblock \href {https://doi.org/10.1016/j.cma.2023.116230}
  {\path{doi:10.1016/j.cma.2023.116230}}.

\bibitem{RN29}
Y.~Liu, F.~Han, L.~Zhang, An extended fictitious node method for surface effect
  correction of bond-based peridynamics, Engineering Analysis with Boundary
  Elements 143 (2022) 78--94.
\newblock \href {https://doi.org/10.1016/j.enganabound.2022.05.023}
  {\path{doi:10.1016/j.enganabound.2022.05.023}}.

\bibitem{RN30}
J.~Zhang, F.~Han, Z.~Yang, J.~Cui, Coupling of an atomistic model and
  bond-based peridynamic model using an extended arlequin framework, Computer
  Methods in Applied Mechanics and Engineering 403 (2023) 115663.
\newblock \href {https://doi.org/10.1016/j.cma.2022.115663}
  {\path{doi:10.1016/j.cma.2022.115663}}.

\bibitem{RN31}
W.~K. Sun, B.~B. Yin, A.~Akbar, V.~K.~R. Kodur, K.~M. Liew, A coupled 3d
  thermo-mechanical peridynamic model for cracking analysis of homogeneous and
  heterogeneous materials, Computer Methods in Applied Mechanics and
  Engineering 418 (2024) 116577.
\newblock \href {https://doi.org/10.1016/j.cma.2023.116577}
  {\path{doi:10.1016/j.cma.2023.116577}}.

\bibitem{RN32}
Y.~Z. Huang, Y.~Yu, Y.~L. Hu, Z.~Y. Yao, D.~Wu, Multi-physical analysis of
  ablation for c/c composites based on peridynamics, Acta Astronautica 214
  (2024) 1--10.
\newblock \href {https://doi.org/10.1016/j.actaastro.2023.10.017}
  {\path{doi:10.1016/j.actaastro.2023.10.017}}.

\bibitem{RN33}
L.~Zhou, Z.~Zhu, A coupled thermo-mechanical peridynamic model for fracture
  behavior of granite subjected to heating and water-cooling processes, Journal
  of Rock Mechanics and Geotechnical Engineering 16~(6) (2024) 2006--2018.
\newblock \href {https://doi.org/10.1016/j.jrmge.2023.07.021}
  {\path{doi:10.1016/j.jrmge.2023.07.021}}.

\bibitem{RN56}
X.~Gu, X.~Li, X.~Xia, E.~Madenci, Q.~Zhang, A robust peridynamic computational
  framework for predicting mechanical properties of porous quasi-brittle
  materials, Composite Structures 303 (2023) 116245.
\newblock \href {https://doi.org/10.1016/j.compstruct.2022.116245}
  {\path{doi:10.1016/j.compstruct.2022.116245}}.

\bibitem{RN34}
B.~Kilic, E.~Madenci, An adaptive dynamic relaxation method for quasi-static
  simulations using the peridynamic theory, Theoretical and Applied Fracture
  Mechanics 53~(3) (2010) 194--204.
\newblock \href {https://doi.org/10.1016/j.tafmec.2010.08.001}
  {\path{doi:10.1016/j.tafmec.2010.08.001}}.

\bibitem{RN35}
F.~Han, Z.~Li, A peridynamics-based finite element method (perifem) for
  quasi-static fracture analysis, Acta Mechanica Solida Sinica 35~(3) (2022)
  446--460.
\newblock \href {https://doi.org/10.1007/s10338-021-00307-y}
  {\path{doi:10.1007/s10338-021-00307-y}}.

\bibitem{RN36}
H.~Wang, L.~Wu, J.~Guo, C.~Yu, Y.~Li, Y.~Wu, Three-dimensional modeling and
  analysis of anisotropic materials with quasi-static deformation and dynamic
  fracture in non-ordinary state-based peridynamics, Applied Mathematical
  Modelling 125 (2024) 625--648.
\newblock \href {https://doi.org/10.1016/j.apm.2023.09.016}
  {\path{doi:10.1016/j.apm.2023.09.016}}.

\bibitem{RN37}
D.~Huang, G.~Lu, P.~Qiao, An improved peridynamic approach for quasi-static
  elastic deformation and brittle fracture analysis, International Journal of
  Mechanical Sciences 94-95 (2015) 111--122.
\newblock \href {https://doi.org/10.1016/j.ijmecsci.2015.02.018}
  {\path{doi:10.1016/j.ijmecsci.2015.02.018}}.

\bibitem{RN38}
J.~Zhong, F.~Han, L.~Zhang, Accelerated peridynamic computation on gpu for
  quasi-static fracture simulations, Journal of Peridynamics and Nonlocal
  Modeling 6~(1) (2023) 206--229.
\newblock \href {https://doi.org/10.1007/s42102-023-00095-8}
  {\path{doi:10.1007/s42102-023-00095-8}}.

\bibitem{RN39}
Q.~Ma, D.~Huang, L.~Wu, D.~Chen, An improved peridynamic model for quasi-static
  and dynamic fracture and failure of reinforced concrete, Engineering Fracture
  Mechanics 289 (2023) 109459.
\newblock \href {https://doi.org/10.1016/j.engfracmech.2023.109459}
  {\path{doi:10.1016/j.engfracmech.2023.109459}}.

\bibitem{RN40}
L.~F. Friedrich, I.~Iturrioz, A.~B. Colpo, S.~Vantadori, Fracture failure of
  quasi-brittle materials by a novel peridynamic model, Composite Structures
  323 (2023) 117402.
\newblock \href {https://doi.org/10.1016/j.compstruct.2023.117402}
  {\path{doi:10.1016/j.compstruct.2023.117402}}.

\bibitem{RN41}
Z.~Chen, X.~Chu, Q.~Duan, A micromorphic peridynamic model and the fracture
  simulations of quasi-brittle material, Engineering Fracture Mechanics 271
  (2022) 108631.
\newblock \href {https://doi.org/10.1016/j.engfracmech.2022.108631}
  {\path{doi:10.1016/j.engfracmech.2022.108631}}.

\bibitem{RN42}
S.~N. Butt, G.~Meschke, Peridynamic analysis of dynamic fracture: influence of
  peridynamic horizon, dimensionality and specimen size, Computational
  Mechanics 67~(6) (2021) 1719--1745.
\newblock \href {https://doi.org/10.1007/s00466-021-02017-1}
  {\path{doi:10.1007/s00466-021-02017-1}}.

\bibitem{RN43}
J.~T. Foster, S.~A. Silling, W.~W. Chen, State based peridynamic modeling of
  dynamic fracture, in: SEM Annual Conf and Exposition on Experimental and
  Applied Mechanics, Albuquerque, USA, 2009, pp. 2312--2317.

\bibitem{RN44}
F.~Bobaru, Influence of van der waals forces on increasing the strength and
  toughness in dynamic fracture of nanofibre networks: a peridynamic approach,
  Modelling and Simulation in Materials Science and Engineering 15~(5) (2007)
  397--417.
\newblock \href {https://doi.org/10.1088/0965-0393/15/5/002}
  {\path{doi:10.1088/0965-0393/15/5/002}}.

\bibitem{RN45}
J.~Yang, Z.~Shen, J.~Zhang, L.~Zhao, An improved implicit non-ordinary
  state-based peridynamics model for modeling crack propagation and coalescence
  problems, Engineering Fracture Mechanics 292 (2023) 109640.
\newblock \href {https://doi.org/10.1016/j.engfracmech.2023.109640}
  {\path{doi:10.1016/j.engfracmech.2023.109640}}.

\bibitem{RN46}
N.~A. Hashim, W.~M. Coombs, C.~E. Augarde, G.~Hattori, An implicit non-ordinary
  state-based peridynamics with stabilised correspondence material model for
  finite deformation analysis, Computer Methods in Applied Mechanics and
  Engineering 371 (2020) 113304.
\newblock \href {https://doi.org/10.1016/j.cma.2020.113304}
  {\path{doi:10.1016/j.cma.2020.113304}}.

\bibitem{RN47}
P.~Li, Z.~Hao, S.~Yu, W.~Zhen, Implicit implementation of the stabilized
  non‐ordinary state‐based peridynamic model, International Journal for
  Numerical Methods in Engineering 121~(4) (2019) 571--587.
\newblock \href {https://doi.org/10.1002/nme.6234}
  {\path{doi:10.1002/nme.6234}}.

\bibitem{RN48}
T.~Ni, M.~Zaccariotto, Q.-Z. Zhu, U.~Galvanetto, Static solution of crack
  propagation problems in peridynamics, Computer Methods in Applied Mechanics
  and Engineering 346 (2019) 126--151.
\newblock \href {https://doi.org/10.1016/j.cma.2018.11.028}
  {\path{doi:10.1016/j.cma.2018.11.028}}.

\bibitem{RN49}
M.~Zaccariotto, F.~Luongo, G.~sarego, U.~Galvanetto, Examples of applications
  of the peridynamic theory to the solution of static equilibrium problems, The
  Aeronautical Journal 119~(1216) (2016) 677--700.
\newblock \href {https://doi.org/10.1017/s0001924000010770}
  {\path{doi:10.1017/s0001924000010770}}.

\bibitem{RN50}
M.~D. Brothers, J.~T. Foster, H.~R. Millwater, A comparison of different
  methods for calculating tangent-stiffness matrices in a massively parallel
  computational peridynamics code, Computer Methods in Applied Mechanics and
  Engineering 279 (2014) 247--267.
\newblock \href {https://doi.org/10.1016/j.cma.2014.06.034}
  {\path{doi:10.1016/j.cma.2014.06.034}}.

\bibitem{RN51}
M.~S. Breitenfeld, P.~H. Geubelle, O.~Weckner, S.~A. Silling, Non-ordinary
  state-based peridynamic analysis of stationary crack problems, Computer
  Methods in Applied Mechanics and Engineering 272 (2014) 233--250.
\newblock \href {https://doi.org/10.1016/j.cma.2014.01.002}
  {\path{doi:10.1016/j.cma.2014.01.002}}.

\bibitem{RN52}
R.~W. Macek, S.~A. Silling, Peridynamics via finite element analysis, Finite
  Elements in Analysis and Design 43~(15) (2007) 1169--1178.
\newblock \href {https://doi.org/10.1016/j.finel.2007.08.012}
  {\path{doi:10.1016/j.finel.2007.08.012}}.

\bibitem{RN53}
J.~Bezanson, A.~Edelman, S.~Karpinski, V.~B. J. S.~r. Shah, Julia: A fresh
  approach to numerical computing, SIAM Review 59~(1) (2017) 65--98.
\newblock \href {https://doi.org/10.1137/141000671}
  {\path{doi:10.1137/141000671}}.

\bibitem{RN54}
U.~Altay, M.~Dorduncu, S.~Kadioglu, An improved peridynamic approach for
  fatigue analysis of two dimensional functionally graded materials,
  Theoretical and Applied Fracture Mechanics 128 (2023) 104152.
\newblock \href {https://doi.org/10.1016/j.tafmec.2023.104152}
  {\path{doi:10.1016/j.tafmec.2023.104152}}.

\bibitem{RN55}
D.~Dipasquale, G.~Sarego, M.~Zaccariotto, U.~Galvanetto, Dependence of crack
  paths on the orientation of regular 2d peridynamic grids, Engineering
  Fracture Mechanics 160 (2016) 248--263.
\newblock \href {https://doi.org/10.1016/j.engfracmech.2016.03.022}
  {\path{doi:10.1016/j.engfracmech.2016.03.022}}.

\end{thebibliography}
\end{document}
\endinput